\documentclass[12pt]{amsart}
\usepackage{amscd,amsmath,amsthm,amssymb}
\usepackage{mathtools}
\usepackage[margin=2cm]{geometry}
\usepackage{epsfig}
\usepackage{rawfonts}
\usepackage{enumerate}
\usepackage{graphics}
\usepackage{multirow}
\usepackage{xspace}
\usepackage[ruled,vlined]{algorithm2e}
\usepackage{graphicx}
\usepackage{mathrsfs}
\usepackage{amsmath}
\usepackage{amsfonts}
\usepackage{amssymb}
\usepackage{amsthm}
\usepackage{graphicx}
\usepackage{booktabs}
\usepackage{caption}
\usepackage{listings}
\usepackage{setspace}
\usepackage[mathscr]{eucal}
\usepackage{pgfplots}
\usepackage{hyperref}
\usepackage{wrapfig}
\usepackage{floatflt,epsfig}
\usepackage{ dsfont }
\usepackage{amscd}
\usepackage{tikz-cd}
\usepackage{fancyhdr}
\usepackage[all]{xy}
\usepackage{latexsym}
\usepackage{amscd}
\usepackage{pifont}
\usepackage{subfig}
\usepackage{easyReview}
\usepackage{subfig}
\usepackage{pstricks-add}
\usepackage{pgf,tikz,pgfplots}
\pgfplotsset{compat=1.15}
\usepackage{mathrsfs}
\usetikzlibrary{arrows}
\usetikzlibrary[patterns]
\usepackage{listings}
\usepackage{verbatim}

\usepackage[normalem]{ulem}
 
\def\NZQ{\Bbb}               % the font for N,Z,Q,R,C
\def\NN{{\NZQ N}}

\newcommand{\cA}{\mathcal{A}}

\newcommand{\cC}{\mathcal{C}}

\newcommand{\cP}{\mathcal{P}}
\newcommand{\Pc}{\mathcal{P}}
\newcommand{\Qc}{\mathcal{Q}}

\newcommand{\cF}{\mathcal{F}}

\newcommand{\cW}{\mathcal{W}}
\newcommand{\cR}{\mathcal{R}}

\newcommand{\cS}{\mathcal{S}}

\renewcommand{\qedsymbol}{$\square$}
\newcommand{\rev}{\mathrm{rev}}

\def\opn#1#2{\def#1{\operatorname{#2}}} % to make operators
\opn\chara{char} \opn\length{\ell} \opn\pd{pd} \opn\rk{rk}
\opn\projdim{proj\,dim} \opn\injdim{inj\,dim} \opn\rank{rank}
\opn\depth{depth} \opn\grade{grade} \opn\height{height}
\opn\embdim{emb\,dim} \opn\codim{codim}

\opn\Tr{Tr} \opn\bigrank{big\,rank}
\opn\superheight{superheight}\opn\lcm{lcm}
\opn\trdeg{tr\,deg}%\emph{
	\opn\reg{reg} \opn\lreg{lreg} \opn\ini{in} \opn\lpd{lpd}
	\opn\size{size} \opn\sdepth{sdepth}
	\opn\link{link}\opn\fdepth{fdepth}\opn\lex{lex}\opn\dist{dist}
	\opn\div{div} \opn\Div{Div} \opn\cl{cl} \opn\Cl{Cl}
	\opn\Spec{Spec} \opn\Supp{Supp} \opn\supp{supp} \opn\Sing{Sing}
	\opn\Ass{Ass} \opn\Min{Min}\opn\Mon{Mon}
	\opn\Ann{Ann} \opn\Rad{Rad} \opn\Soc{Soc}
	\opn\Im{Im} \opn\Ker{Ker} \opn\Coker{Coker} \opn\Am{Am}
	\opn\Hom{Hom} \opn\Tor{Tor} \opn\Ext{Ext} \opn\End{End}
	\opn\Aut{Aut} \opn\id{id}
	
	\opn\nat{nat}
	\opn\pff{pf}%   \pf exists already
	\opn\Pf{Pf} \opn\GL{GL} \opn\SL{SL} \opn\mod{mod} \opn\ord{ord}
	\opn\Gin{Gin} \opn\Hilb{Hilb}\opn\sort{sort}
	\opn\aff{aff} \opn
	\con{conv} \opn\relint{relint} \opn\st{st}
	\opn\lk{lk} \opn\cn{cn} \opn\core{core} \opn\vol{vol}
	\opn\link{link} \opn\star{star}\opn\lex{lex}\opn\set{set}
	\opn\gr{gr}
	
	\def\pot#1#2{#1[\kern-0.28ex[#2]\kern-0.28ex]}

	\opn\dirlim{\underrightarrow{\lim}}
	\opn\inivlim{\underleftarrow{\lim}}
	\let\to=\rightarrow
	
	\def\Implies{\ifmmode\Longrightarrow \else
		\unskip${}\Longrightarrow{}$\ignorespaces\fi}
	\def\implies{\ifmmode\Rightarrow \else
		\unskip${}\Rightarrow{}$\ignorespaces\fi}
	\def\iff{\ifmmode\Longleftrightarrow \else
		\unskip${}\Longleftrightarrow{}$\ignorespaces\fi}

	\let\:=\colon
	\let\epsilon\varepsilon
	\let\kappa=\varkappa
	\def\qed{\ifhmode\textqed\fi
		\ifmmode\ifinner\quad\qedsymbol\else\dispqed\fi\fi}
	\def\textqed{\unskip\nobreak\penalty50
		\hskip2em\hbox{}\nobreak\hfil\qedsymbol
		\parfillskip=0pt \finalhyphendemerits=0}
	\def\dispqed{\rlap{\qquad\qedsymbol}}
	
	\opn\dis{dis}
	\def\pnt{{\raise0.5mm\hbox{\large\bf.}}}
	
	\opn\Lex{Lex}
	
        \newtheorem{Theorem}{Theorem}[section]
	
	\newtheorem{Corollary}[Theorem]{Corollary}
	\newtheorem{Proposition}[Theorem]{Proposition}
	\newtheorem{Remark}[Theorem]{Remark}

	\newtheorem{Definition}[Theorem]{Definition}
	
	\newtheorem{Conjecture}[Theorem]{Conjecture}
	\newtheorem{Question}[Theorem]{Question}

\begin{document}

    \title[Advances in Polyomino Ideals]{Recent Advances in the Theory of Polyomino Ideals}

    \author{Francesco Navarra}
    \address{Sabanci University, Faculty of Engineering and Natural Sciences, Orta Mahalle, Tuzla 34956, Istanbul, Turkey}
    \email{francesco.navarra@sabanciuniv.edu}
    
    \author{Ayesha Asloob Qureshi}      
    \address{Sabanci University, Faculty of Engineering and Natural Sciences, Orta Mahalle, Tuzla 34956, Istanbul, Turkey}	
    \email{aqureshi@sabanciuniv.edu, ayesha.asloob@sabanciuniv.edu}

\thanks{The authors are supported by Scientific and Technological Research Council of Turkey T\"UB\.{I}TAK under the Grant No: 124F113, and are thankful to T\"UB\.{I}TAK for their support.The first author is member of INDAM-GNSAGA and acknowledge its support. }



\maketitle
\begin{abstract}
Polyomino ideals, defined as the ideals generated by the inner $2$-minors of a polyomino, are a class of binomial ideals whose algebraic properties are closely related to the combinatorial structure of the underlying polyomino. We provide a unified account of recent advances on two central themes: the characterization of prime polyomino ideals and the emerging connection between the Hilbert–Poincaré series and Gorensteinness of \(K[\Pc]\) with the classical rook theory. Some further related properties, as radicality, primary decomposition, and levelness are discussed, and a \texttt{Macaulay2} package, namely \texttt{PolyominoIdeals}, is also presented.
\end{abstract}

\section*{Introduction}
Let $X = (x_{ij})_{i=1,\dots,m\atop j=1,\dots,n}$ be an $m \times n$ matrix of
indeterminates over a field $K$.  
The study of the ideal generated by all $t$--minors of $X$ is a classical topic 
in Commutative Algebra and Algebraic Geometry.  
For any integer $1 \le t \le \min\{m,n\}$, the ideal generated by all $t$--minors, 
called a \emph{determinantal ideal}, has been extensively investigated.  
These ideals play a fundamental role in Algebraic Geometry, as they define 
several classical varieties, including the Veronese and Segre embeddings.  
A standard reference on determinantal ideals and their associated algebras is 
the monograph of Bruns and Vetter~\cite{BV}.

Over time, various generalizations of determinantal ideals have been introduced, 
including one-sided and two-sided ladder determinantal ideals.  
However, when $I$ is an ideal generated by an arbitrary set of $t$--minors of $X$, 
the situation becomes far more intricate, even in the smallest nontrivial case 
$t=2$.  
A central problem is to determine when such an ideal is \emph{prime} or 
\emph{radical}, and to describe its primary decomposition.  
It is shown in~\cite[Corollary~2.2]{HHHKR} that $I$ is always radical when 
$X$ is a $2 \times n$ matrix, and in this case the authors also provide an 
explicit minimal primary decomposition.  
The situation becomes considerably more subtle as soon as either $m$ or $n$ is 
at least $3$: explicit examples show that such ideals need not be radical in 
general.  
This leads naturally to the problem of characterizing those arbitrary sets of 
$2$--minors whose ideal is radical.  
Its significance in Algebraic Statistics (see~\cite{HHHKR}) has motivated a 
systematic investigation of ideals generated by arbitrary subsets of $2$--minors 
of an $m \times n$ matrix.

Within this broader framework, \emph{polyomino ideals} arise as a particularly 
rich and structured subclass of ideals of 2-minors. A polyomino is a finite union of unit squares in the plane joined edge to edge. Their enumerative and structural properties link them to tilings, lattice-path combinatorics, and discrete geometry (see \cite{Golomb1994}). The systematic study of polyominoes, and more generally of collections of cells, from the perspective of Commutative Algebra, was initiated by the second author in \cite{Qureshi2012}. To each collection of cells $\Pc$ one associates
the binomial ideal $I_{\Pc}$ generated by the inner $2$-minors of $\cP$ in the polynomial ring $S_{\Pc}=K[x_a : a \in V(\Pc)]$, where $V(\Pc)$ is the vertex set of $\Pc$ and $K$ is a field. The corresponding quotient $K[\Pc]=S_{\Pc}/I_{\Pc}$ is called the coordinate
ring of $\Pc$.  This construction fits naturally into the general framework of binomial
ideals and affine semigroup rings (see, for example, \cite{HHO_Book}), and leads to the
central theme of relating algebraic properties of $K[\Pc]$ to combinatorial features of
$\Pc$.

The goal of this paper is to present a unified account of the current understanding of the primality, radicality, and Hilbert–Poincaré series of polyomino ideals.
We review the characterization of prime polyomino ideals, recent progress on radicality and primary decomposition, and the surprising appearance of classical combinatorial invariants, rook and switching rook polynomials, in the description of Hilbert–Poincaré series and Gorensteinness of coordinate rings of polyominoes.
Special attention is devoted to the combinatorial structures driving these phenomena and to open problems at the forefront of current research.

A breakdown of this paper is as follows.
Section \ref{Section: Preliminaries} collects the basic terminology on
collections of cells, inner 2-minors, and the construction of the polyomino ideal \(I_{\Pc}\).
Sections \ref{Section:SimplePolyominoes} and \ref{sec:non-simple} are devoted to
the study of primality of polyomino ideals.
Whether \(I_{\Pc}\) is prime is strongly influenced by the topological structure of the underlying polyomino, most notably, by the presence or absence of holes.
Polyominoes that contain no holes are called simple polyominoes, and their primality theory is treated in Section \ref{Section:SimplePolyominoes}. We also discuss the toric representation of their coordinate rings via edge rings of bipartite graphs, their resulting normality and Koszulness.

In Section~\ref{sec:non-simple} we turn to non-simple polyominoes.
We first discuss the obstruction to obtaining a toric representation in this
setting, as conjectured in \cite{HibiQureshi2015}, and the localization
techniques used in that work.  
We then outline the toric description of Shikama~\cite{Shikama2015} and its role
in subsequent developments.  
A central theme of this section is the combinatorial criterion introduced by
Mascia, Rinaldo and Romeo~\cite{MasciaRinaldoRomeo2020}, based on zig--zag
walks, which could provide a powerful tool for characterizing non-prime polyominoes.
We present the complete classification obtained for closed path and weakly
closed path polyominoes, together with related results for grid 
polyominoes \cite{CistoNavarra2021,CNU1,CNU2,CNU3,DinuNavarra2025,Navarra2025}.

In Section \ref{Section:Radicality}, we review what is currently known about the radicality and primary decomposition of polyomino ideals, topics that remain far less developed than primeness.
We highlight the few available positive results (for example, radicality of certain cross polyominoes) as well as the polyocollection framework, which provides the first general methodology for describing minimal primary decompositions of non-prime polyomino ideals.

In Section \ref{Section:RookPoly}, we focus on Hilbert–Poincaré series and their unexpectedly rich connection with rook theory.
Rook polynomials were first introduced in classical enumerative combinatorics to count non-attacking rook placements on a chessboard.
They have since appeared in a wide range of areas: permutation enumeration, inclusion-exclusion, statistical mechanics, and the study of Ferrers boards and permutation statistics.
What makes them particularly significant in the context of polyomino ideals is that they encode the $h$-polynomial of the coordinate ring of a collection of cells.
We describe the connection, first established for $L$-convex polyominoes in
\cite{EneHerzogQureshiRomeo2021}, between the Castelnuovo-Mumford regularity of $K[\Pc]$ and the
rook number of $\Pc$, and the corresponding relation between the
$h$-polynomial of $K[\cP]$ and the rook polynomial of $\cP$.  
We then discuss the extension of these ideas to simple thin polyominoes,
closed paths, weakly closed paths, and grid polyominoes, as well as the introduction of the switching rook polynomial, which refines the rook polynomial in the non-thin case. 

In Section~\ref{Sectio: canonical module}, we survey recent advances concerning the canonical module and the related generalizations of the Gorenstein property, namely pseudo-Gorensteiness and levelness. We first present the characterizations obtained in \cite{RinaldoRomeoSarkar2024} for simple paths. We then discuss the combinatorial description of the canonical module of circle closed paths, developed in \cite{DinuNavarra2025}. Finally, we consider a particular family of Ferrers diagrams whose Cohen–Macaulay type is provided by Fuss–Catalan numbers, as established in \cite{AlinStefan}. 

We conclude with Section~\ref{Section: Package}, where we illustrate the combinatorial and algebraic functions implemented in the \texttt{Macaulay2} package \texttt{PolyominoIdeals} \cite{Package_M2, M2}.

%%%%%%%%%%%%%%%%%%%%%%%%%%%%%%%%%%%%%%%%%%%%%%
%%%%%%%%%%%%%%%%%%%%%%%%%%%%%%%%%%%%%%%%%%%%%%%%%%%%%%%%%%%%%%%%%%%%%%%%%%%%%%%%%%%%%%%%%%%%%%%%%%%%%%%%%%%%%%%%%%%%%%%%%%%%%%%%%%%%%%%%%%%%%%%%%%%%%%%%%%%%%%%%%%%%%%%%%%%%%%%%%%%%%%%%%%%%%%%%%%%%%%%%%%%%%%%%%%%%%%%%%%%%%%%%%%%%%%%%%%%%%%%%%%%%%%%%%%%%%%%%%%%%%%%%%%%%%%%%%%%%%%%%%%%%%%%%%%%%%%%%%%%%%%%%%%%%%%%%%%%%%%%%%%

\section{Required terminologies related to collection of cells}\label{Section: Preliminaries}

This section is devoted to introducing the definitions and notations related to collections of cells and the associated ideals of $2$-minors.

Let $(i,j), (k,l) \in \mathbb{Z}^2$. We define a partial order on $\mathbb{Z}^2$ by setting $(i,j) \leq (k,l)$ if and only if $i \leq k$ and $j \leq l$. Given $a = (i,j)$ and $b = (k,l)$ in $\mathbb{Z}^2$ with $a \leq b$, we define the set
\[
[a,b] = \{(m,n) \in \mathbb{Z}^2 \mid i \leq m \leq k,\ j \leq n \leq l\}
\]
as an \emph{interval} in $\mathbb{Z}^2$. If $i < k$ and $j < l$, then $[a,b]$ is called a \emph{proper} interval. In this case, we refer to $a$ and $b$ as the \emph{diagonal corners} of $[a,b]$, and define $c = (i,l)$ and $d = (k,j)$ as the \emph{anti-diagonal corners}. If $j = l$ (respectively, $i = k$), then $a$ and $b$ are said to be in a \emph{horizontal} (respectively, \emph{vertical}) position.

A proper interval $C = [a,b]$ with $b = a + (1,1)$ is called a \emph{cell} of $\mathbb{Z}^2$. The points $a$, $b$, $c$, and $d$ are referred to as the \emph{lower-left}, \emph{upper-right}, \emph{upper-left}, and \emph{lower-right} corners of $C$, respectively. We denote the set of \emph{vertices} and \emph{edges} of $C$ by $V(C) = \{a, b, c, d\}, \quad E(C) = \{\{a,c\}, \{c,b\}, \{b,d\}, \{a,d\}\}.$ For a collection of cells $\cP$ in $\mathbb{Z}^2$, the sets of vertices and edges are defined as $V(\cP) = \bigcup_{C \in \cP} V(C), \quad E(\cP) = \bigcup_{C \in \cP} E(C).$ The \emph{rank} of $\cP$, denoted by $|\cP|$, is the number of cells in $\cP$. By convention, the empty set is considered a collection of cells of rank $0$.

Consider two cells $A$ and $B$ in $\mathbb{Z}^2$, with lower-left corners $a = (i,j)$ and $b = (k,l)$, respectively, and suppose that $a \leq b$. The \emph{cell interval} $[A,B]$, also referred to as a \emph{rectangle}, is the set of all cells in $\mathbb{Z}^2$ whose lower-left corners $(r,s)$ satisfy $i \leq r \leq k$ and $j \leq s \leq l$.

\begin{figure}[h]
    \centering
    \subfloat[]{\includegraphics[scale=0.4]{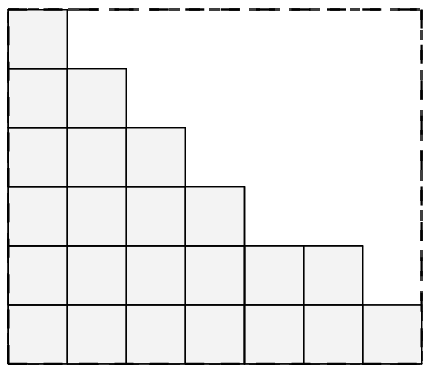}}\quad
    \subfloat[]{\includegraphics[scale=0.4]{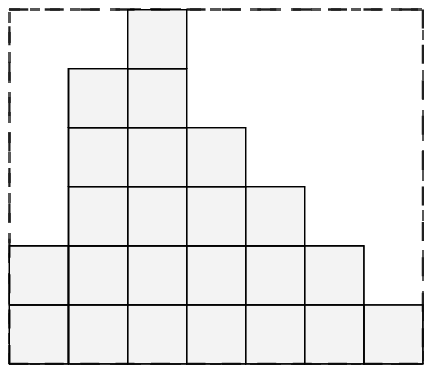}}\quad
    \subfloat[]{\includegraphics[scale=0.4]{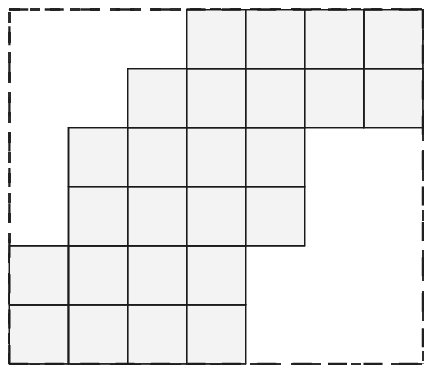}}\quad
    \subfloat[]{\includegraphics[scale=0.4]{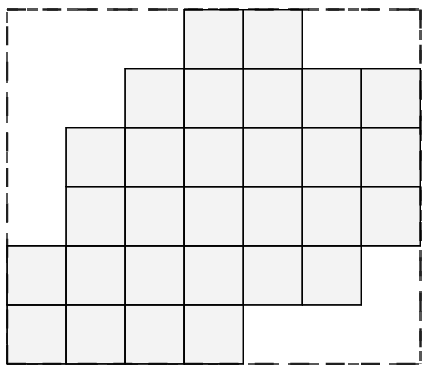}}
    \caption{A Ferrer diagram, a stack, a parallelogram and a directed convex.}
    \label{fig:Ferrer diagram and stack}
\end{figure}

Let $\Pc$ be a collection of cells. The cell interval $[A,B]$ is called the \emph{minimal bounding rectangle} of $\Pc$ if $\Pc \subseteq [A,B]$ and there exists no other rectangle in $\mathbb{Z}^2$ that properly contains $\Pc$ and is properly contained in $[A,B]$ (see Figure~\ref{fig:Ferrer diagram and stack}). If the corners $(i,j)$ and $(k,l)$ are in horizontal (respectively, vertical) position, we say that the cells $A$ and $B$ are in \emph{horizontal} (respectively, \emph{vertical}) position. 

An interval $[a,b]$ with $a=(i,j)$, $b=(k,j)$, and $i<k$ is called a \textit{horizontal edge interval} of $\mathcal{P}$ if the sets $\{(\ell,j),(\ell+1,j)\}$ are edges of cells of $\mathcal{P}$ for all $\ell=i,\dots,k-1$. 
If $\{(i-1,j),(i,j)\}$ and $\{(k,j),(k+1,j)\}$ do not belong to $E(\mathcal{P})$, then $[a,b]$ is a \textit{maximal horizontal edge interval} of $\mathcal{P}$. 
Vertical and maximal vertical edge intervals are defined analogously.%The cell interval $[A,B]$ is called a \textit{column} (respectively, a \textit{row}) of $\cP$ if the cells $A$ and $B$ are in vertical (respectively, horizontal) position, all the cells in $[A,B]$ belong to $\cP$, and there is no cell interval $[A',B']$ with $A'$ and $B'$ in vertical (respectively, horizontal) position such that $[A,B]\subsetneq [A',B']$. The collection of cells in Figure~\ref{Figure: Polyomino + weakly conn. coll. of cells + three disconnected comp.} (A) has thirteen columns and ten rows.

\begin{figure}[h]
    \centering
    \subfloat[]{\includegraphics[scale=0.5]{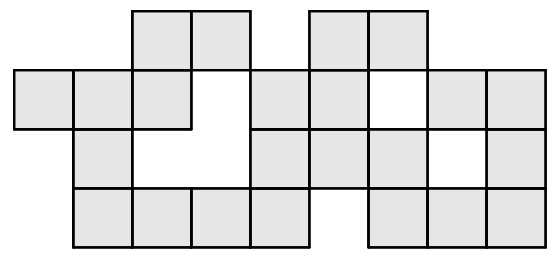}}
    \subfloat[]{\includegraphics[scale=0.5]{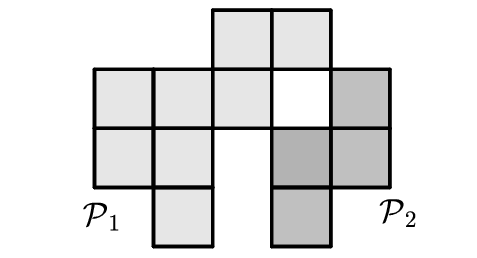}}
    \subfloat[]{\includegraphics[scale=0.5]{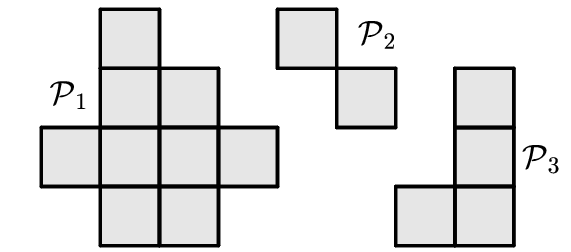}}
    \caption{A polyomino, a weakly connected collection of cells with two connected components, and a collection of cells with three weakly connected components.}
    \label{Figure: Polyomino + weakly conn. coll. of cells + three disconnected comp.}
\end{figure}

A finite collection of cells $\Pc$ is said to be \emph{weakly connected} if, for any two cells $C$ and $D$ in $\Pc$, there exists a sequence of cells $\mathcal{C} \colon C = C_1, \dots, C_m = D$ in $\Pc$ such that $V(C_i) \cap V(C_{i+1}) \neq \emptyset$ for all $i = 1, \dots, m-1$. For an illustration, see Figure~\ref{Figure: Polyomino + weakly conn. coll. of cells + three disconnected comp.} (B).

If $\Pc = \bigcup_{i=1}^s \Pc_i$, where each $\Pc_i$ is a weakly connected collection of cells and $V(\Pc_i) \cap V(\Pc_j) = \emptyset$ for all $i \neq j$, then $\Pc_1, \dots, \Pc_s$ are called the \emph{weakly connected components} of $\Pc$. The collection of cells in Figure~\ref{Figure: Polyomino + weakly conn. coll. of cells + three disconnected comp.} (C) has three weakly connected components.

A finite collection of cells $\Pc$ is called \emph{connected}, or simply a \emph{polyomino}, if for any two cells $C$ and $D$ in $\Pc$, there exists a sequence of cells $\mathcal{C} \colon C = C_1, \dots, C_m = D$ in $\Pc$ such that $C_i \cap C_{i+1}$ is an edge shared by both $C_i$ and $C_{i+1}$ for all $i = 1, \dots, m-1$. Such a sequence is called a \emph{path} from $C$ to $D$ in $\Pc$. An example of a polyomino is shown in Figure~\ref{Figure: Polyomino + weakly conn. coll. of cells + three disconnected comp.} (A). Moreover, if we denote by $(a_i,b_i)$ the lower left corner of $C_i$ for all $i=1,\dots,m$, then $\cC$ has a \textit{change of direction} at $C_k$, for some $2\leq k \leq m-1$, if $a_{k-1} \neq a_{k+1}$ and $b_{k-1} \neq b_{k+1}$; in this case, $\{C_{k-1},C_k,C_{k+1}\}$ is said to be the set of the \textit{cells of a change of direction}. 

A subcollection $\Pc' \subseteq \Pc$ is called a \emph{connected component} of $\Pc$ if $\Pc'$ is a polyomino and is maximal with respect to set inclusion; that is, for any $A \in \Pc \setminus \Pc'$, the union $\Pc' \cup \{A\}$ is not a polyomino. For instance, the collection of cells in Figure~\ref{Figure: Polyomino + weakly conn. coll. of cells + three disconnected comp.} (B) has two connected components $\Pc_1$ and $\Pc_2$.

A collection of cells $\Pc$ is said to be \emph{simple} if, for any two cells $C$ and $D$ not in $\Pc$, there exists a path of cells in the complement of $\Pc$ connecting $C$ to $D$. Roughly speaking, simple polyominoes are the polyominoes without holes. The collections in Figures~\ref{Figure: Polyomino + weakly conn. coll. of cells + three disconnected comp.} (A) and (B) are not simple, while the one in (C) is.

A collection of cells $\Pc$ is said to be \emph{row convex} (respectively, \emph{column convex}) if, for any two cells $A$ and $B$ of $\Pc$ in horizontal (respectively, vertical) alignment, the entire cell interval $[A, B]$ lies in $\Pc$. If $\Pc$ is both row and column convex, it is called \emph{convex}. In Figure~\ref{Figure: Polyomino + weakly conn. coll. of cells + three disconnected comp.} (C), the collection $\Pc = \Pc_1 \cup \Pc_2 \cup \Pc_3$ is column convex but not row convex. Moreover, each weakly connected component of $\Pc$ is convex.

Among the convex polyominoes, we have some very well-studied sub-classes. Let $\Pc$ be a convex polyomino with minimal bounding rectangle $[A,B]$. Then:
\begin{enumerate}
    \item $\Pc$ is called a \emph{Ferrer diagram} if at least three corner cells of $[A,B]$ are in $\Pc$ (Figure~\ref{fig:Ferrer diagram and stack} (A)).
    \item $\Pc$ is called a \emph{stack} if two adjacent corner cells of $[A,B]$ belong to $\Pc$ (Figure~\ref{fig:Ferrer diagram and stack} (B)).
    \item $\Pc$ is called a \emph{parallelogram} if two opposite corner cells of $[A,B]$ belong to $\Pc$ (Figure~\ref{fig:Ferrer diagram and stack} (C)).
    \item $\Pc$ is called \emph{directed convex} if at least one corner cell of $[A,B]$ belongs to $\Pc$ (Figure~\ref{fig:Ferrer diagram and stack} (D)).
\end{enumerate}

It follows directly from the above definitions that every Ferrer diagram, stack polyomino, and parallelogram polyomino is a directed convex polyomino.

\subsection{Ideals of inner 2-minors of collections of cells and polyomino ideals.} Let $\Pc$ be a collection of cells, and let $S_\Pc = K[x_v : v \in V(\Pc)]$ be the polynomial ring associated to $\Pc$, where $K$ is a field. A proper interval $[a,b]$ is called an \emph{inner interval} of $\Pc$ if all cells in $\Pc_{[a,b]}$ belong to $\Pc$. If $[a,b]$ is an inner interval of $\Pc$, with diagonal corners $a$ and $b$ and anti-diagonal corners $c$ and $d$, then the binomial $x_a x_b - x_c x_d$ is called an \emph{inner 2-minor} of $\Pc$. The ideal $I_{\Pc} \subset S_\Pc$ generated by all inner 2-minors of $\Pc$ is called the \emph{ideal of 2-minors} of $\Pc$. If $\Pc$ is a polyomino, then $I_{\Pc}$ is referred to as the \emph{polyomino ideal} of $\Pc$. The quotient ring $K[\Pc] = S_\Pc / I_{\Pc}$ is called the \emph{coordinate ring} of $\Pc$.

\begin{figure}[h]
		\centering
		\includegraphics[scale=0.9]{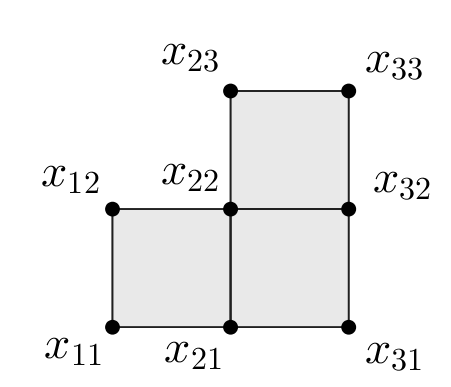}
		\caption{A polyomino $\cP$.}
		\label{Figure: Polyomino introduction}
	\end{figure}

For example, let $\cP$ be the polyomino represented in Figure \ref{Figure: Polyomino introduction}.
	  The ideal $I_\cP$ is generated by the following binomials:
			\[
				x_{11}x_{22} - x_{12}x_{21},\ x_{21}x_{32} - x_{22}x_{31},\ x_{22}x_{33} - x_{23}x_{32},\ x_{11}x_{32} - x_{12}x_{31},\
				x_{21}x_{33} - x_{23}x_{31}.
			\]
	
For ease of notation, we say that a polyomino $\Pc$ enjoys a property $(P)$ if its polyomino ideal $I_{\Pc}$, or equivalently its coordinate ring $K[\Pc]$, satisfies property $(P)$.

\begin{Remark}\label{rem:p=q}
Let $\Pc$ be a collection of cells, and let $\Qc$ be the collection obtained from $\Pc$ by a symmetry of the plane—i.e., by a translation, rotation, reflection, or glide reflection. Then $I_\Pc$ and $I_\Qc$ define the same ideal up to a relabeling of variables; in particular, $K[\Pc] \cong K[\Qc]$.
\end{Remark}

%We conclude this subsection by introducing some notions and results related to the Hilbert-Poincaré series of a graded $K$-algebra $R/I$. Let $R$ be a graded $K$-algebra and $I$ a homogeneous ideal of $R$. Then $R/I$ inherits a natural grading given by $R/I = \bigoplus_{k \in \mathbb{Z}} (R/I)_k.$ The function $\mathrm{H}_{R/I} \colon \mathbb{Z} \to \mathbb{Z}$ defined by $\mathrm{H}_{R/I}(k) = \dim_K (R/I)_k$ is called the \emph{Hilbert function} of $R/I$. The associated formal power series
%\[\mathrm{HP}_{R/I}(t) = \sum_{k \in \mathbb{Z}} \mathrm{H}_{R/I}(k) t^k\]is called the \emph{Hilbert-Poincaré series} of $R/I$. By the Hilbert–Serre theorem, there exists a unique polynomial $h(t) \in \mathbb{Z}[t]$, called the \emph{$h$-polynomial} of $R/I$, such that $h(1) \neq 0$ and
%\[\mathrm{HP}_{R/I}(t) = \frac{h(t)}{(1 - t)^d},\]where $d$ is the Krull dimension of $R/I$. Moreover, if $R/I$ is Cohen–Macaulay, then the Castelnuovo–Mumford regularity satisfies $\mathrm{reg}(R/I) = \deg h(t)$. If $\Pc$ is a collection of cells, we denote the $h$-polynomial of $K[\Pc]$ by $h_{K[\Pc]}(t)$ throughout this paper.

%%%%%%%%%%%%%%%%%%%%%%%%%%%%%%%%%%%%%%%%%%%%%%%%%%%%%%%%%%%%%%%%%%%%%%%%%%%%%%%%%%%%%%%%%%%%%%%%%%%%%%%%%%%%%%%%%%%%%%%%%%%%%%%%%%%%%%%%%%%%%%%%%%%%%%%%%%%%%%%%%%%%%%%%%%%%%%%%%%%%%%%%%%%%%%%%%%%%%%%%%%%%%%%%%%%%%%%%%%%%%%%%%%%%%%%%%%%%%%%%%%%%%%%%%%%%%%%%%%%%%%%%%%%%%%%%%%%%%%%%%%%%%%%%%%%%%%%%%%%%%%%%%%%%%%%%%%%%%%%%%%%%%%%%%%%%%%%%%%%%%%%%%%%%%%%%%%%%%%%%%%%%%%%%%%%%%%%%%%%%%%%%%%%%%%%%%%%%%%%%%%%%%%%%%%%%%%%%%%%%

%
\section{Simple polyominoes and their toric representation}\label{Section:SimplePolyominoes}

Among all classes of polyominoes, simple polyominoes admit the most complete and 
coherent algebraic description.  
Their ideals of inner $2$-minors coincide with toric ideals arising from 
naturally associated bipartite graphs, a correspondence that allows one to 
deduce important structural properties such as normality, Cohen--Macaulayness, 
and Koszulness.  
In this section we present the toric viewpoint for simple polyominoes and review 
the results that place them among the best understood classes of polyominoes.

\subsection*{Edge rings and toric ideals of bipartite graphs}

Let $G$ be a finite simple graph with vertex set $V(G)=\{x_1,\dots,x_m\}$, identified with variables
in the polynomial ring $K[x_1,\dots,x_m]$.  
For each edge $\{x_i,x_j\}\in E(G)$ we introduce a new variable $t_{ij}$, and we consider the
polynomial ring $T = K[t_{ij} : \{x_i,x_j\} \in E(G)].$ The \emph{edge ring} of $G$ is the $K$-subalgebra of the polynomial ring
$K[x_1,\dots,x_m]$ generated by all quadratic monomials corresponding to edges,
\[
K[G] = K[x_i x_j : \{x_i,x_j\} \in E(G)].
\]
Equivalently, $K[G]$ is the image of the $K$-algebra homomorphism
\[
\psi \colon T \longrightarrow K[x_1,\dots,x_m], \qquad
\psi(t_{ij})= x_i x_j .
\]
The kernel $I_G = \ker(\psi)$ is called the \emph{toric ideal} of $G$; it is a prime ideal
generated by binomials corresponding to combinatorial relations among the edges. When $G$ is bipartite, the structure of $I_G$ is particularly transparent. A cycle $C$ of even length in $G$, written as $
x_{i_1}, x_{i_2}, \dots, x_{i_{2r}}, x_{i_1},$ gives rise to a binomial
\[
f_C = t_{i_1 i_2} t_{i_3 i_4} \cdots t_{i_{2r-1} i_{2r}}
      - t_{i_2 i_3} t_{i_4 i_5} \cdots t_{i_{2r} i_1}.
\]
It is well known that for bipartite graphs, the toric ideal 
$I_G$ is generated by the binomials attached to all even cycles of $G$, for example see \cite[Lemma 5.9]{HHO_Book}.  
Moreover, $I_G$ is generated by quadrics if and only if every even cycle of 
length greater than four has a chord, i.e., if and only if $G$ is \emph{weakly chordal}
\cite{OhsugiHibiKoszul1999}. 

\subsection*{Toric representation of simple polyominoes}

Let $\Pc$ be a polyomino.  
Let $\{V_1,\dots,V_m\}$ be its maximal vertical edge intervals and 
$\{H_1,\dots,H_n\}$ its maximal horizontal edge intervals.  
The associated bipartite graph $G_{\Pc}$ has vertex set
\[
V(G_{\Pc})=\{v_1,\dots,v_m\} \sqcup \{h_1,\dots,h_n\},
\]
where $v_i$ corresponds to $V_i$ and $h_j$ to $H_j$.  
There is an edge $\{v_i,h_j\}\in E(G_{\Pc})$ precisely when 
$V_i\cap H_j$ is a vertex of $\Pc$.  
This defines a $K$-algebra homomorphism
\[
\Phi_{\Pc} \colon S_{\Pc} \longrightarrow K[G_{\Pc}], \qquad \text{ with } \qquad
\Phi_\Pc(x_{a}) = v_i h_j \quad \text{whenever } V_i\cap H_j=\{a\}.
\]
The kernel of $\Phi_\Pc$, that is $\ker(\Phi_\Pc)$, is the toric ideal of the bipartite graph $G_{\Pc}$.

A key feature of simple polyominoes is that the combinatorics of their maximal 
edge intervals forces $G_{\Pc}$ to be weakly chordal, so that its toric ideal is 
generated by quadratic binomials corresponding exactly to the inner $2$-minors of $\Pc$.

\begin{Theorem}
\label{Thm:SimpleToric}\cite[Theorem 2.2]{QSS}
If $\Pc$ is a simple polyomino, then $I_{\Pc} = \ker(\Phi_\Pc)$, and $I_\Pc$ is prime. %Moreover, $I_{\Pc}$ admit a quadratic Gr\"obner basis with respect to a suitable monomial order, and $K[\Pc]$ is Koszul. 
\end{Theorem}

The proof relies on two facts:  
(i) the graph $G_{\Pc}$ associated with a simple polyomino is weakly chordal, and  
(ii) the toric ideal of the edge ring of a weakly chordal bipartite graph admits a 
quadratic Gröbner basis \cite{OhsugiHibiKoszul1999}.  
Since these quadratic binomials correspond exactly to the inner $2$-minors of $\Pc$, 
one obtains $I_{\Pc} = J_{\Pc}$.

Before the toric description was established, 
Herzog, Qureshi, and Shikama \cite{HerzogQureshiShikama2015} 
introduced the class of \emph{balanced} polyominoes and showed that for such 
$\Pc$, the ideal $I_{\Pc}$ is the lattice ideal of a saturated lattice; in particular, 
$I_{\Pc}$ is prime and they computed its universal Gröbner basis.  
Later, Herzog and Saeedi Madani \cite{HerzogSaeedi2014} proved that balanced 
polyominoes coincide exactly with simple polyominoes.  
This provides an alternative route to the primality of $I_{\Pc}$ for simple polyominoes.

Combining the results in 
\cite{HerzogQureshiShikama2015,HerzogSaeedi2014,QSS} 
we obtain the following algebraic properties of simple polyominoes.

\begin{enumerate}
    \item $I_{\Pc}$ is prime, and $\operatorname{ht}(I_{\Pc}) = |\Pc|$.
    \item The universal Gr\"obner basis of $I_\Pc$ consists of squarefree binomials.
    \item $I_{\Pc}$ has a quadratic Gr\"obner basis with respect to a suitable order.
    \item $K[\Pc]$ is a normal Cohen--Macaulay domain.
    \item $K[\Pc]$ is Koszul.
\end{enumerate}

Following the same approach as in \cite{QSS}, 
Cisto, Navarra, and Utano \cite[Theorem.~3.3]{CNU2} 
generalized the result of \cite{QSS} to simple and weakly connected collections of cells.

%%%%%%%%%%%%%%%%%%%%%%%%%%%%%%%%%%%%%%%%%%%%%%%%%%%%%%%%%%%%%%%%%%%%%%%%%%%%%%%%%%%%%%%%%%%%%%%%%%%%%%%%%%%%%%%%%%%%%%%%%%%%%%%%%%%%%%%%%%%%%%%%

\section{Characterization of non-simple prime polyominoes}\label{sec:non-simple}

Polyominoes with one or more holes are called \emph{multiply connected}
polyominoes, or \emph{non-simple} polyominoes.  
We will use the terminology ``non-simple'' throughout this section.  
As explained in Section~\ref{Section:SimplePolyominoes}, the coordinate ring of a simple polyomino admits a toric description via the edge ring of a naturally associated bipartite graph. A natural question is whether such a toric representation can be extended to non-simple polyominoes.  
The following conjecture, formulated by Hibi and Qureshi~\cite{HibiQureshi2015}, suggests that this is not possible.

\begin{Conjecture}[\cite{HibiQureshi2015}]
A polyomino ideal $I_\Pc$ arises as the toric ideal of the edge ring of a
finite simple graph if and only if $\Pc$ is simple.
\end{Conjecture}

\begin{figure}[h]
    \centering
    \subfloat[]{\includegraphics[scale=0.6]{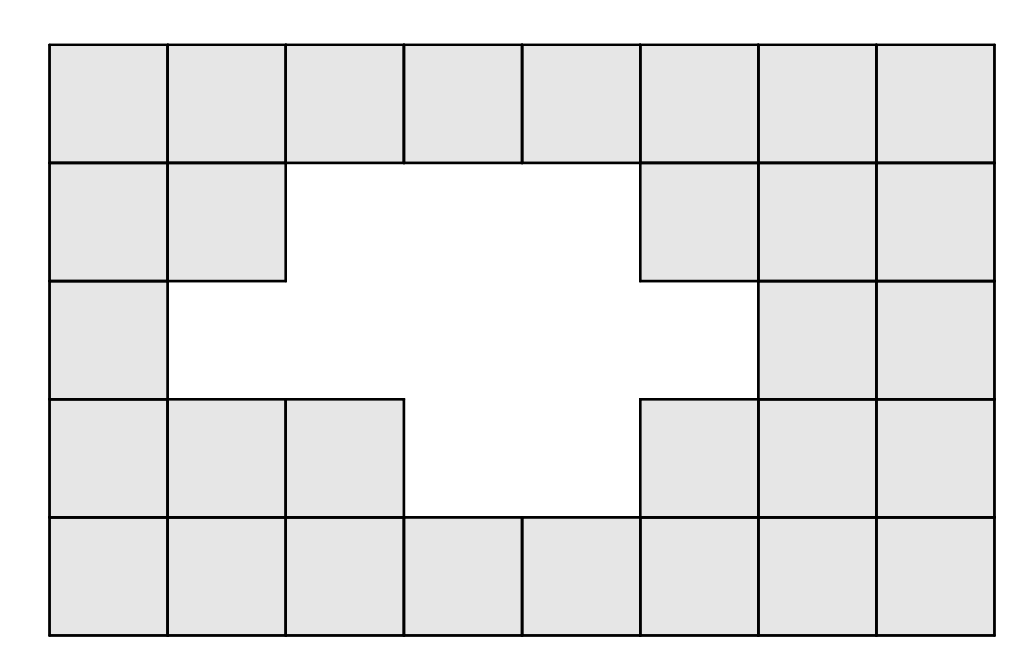}}\qquad
    \subfloat[]{\includegraphics[scale=0.6]{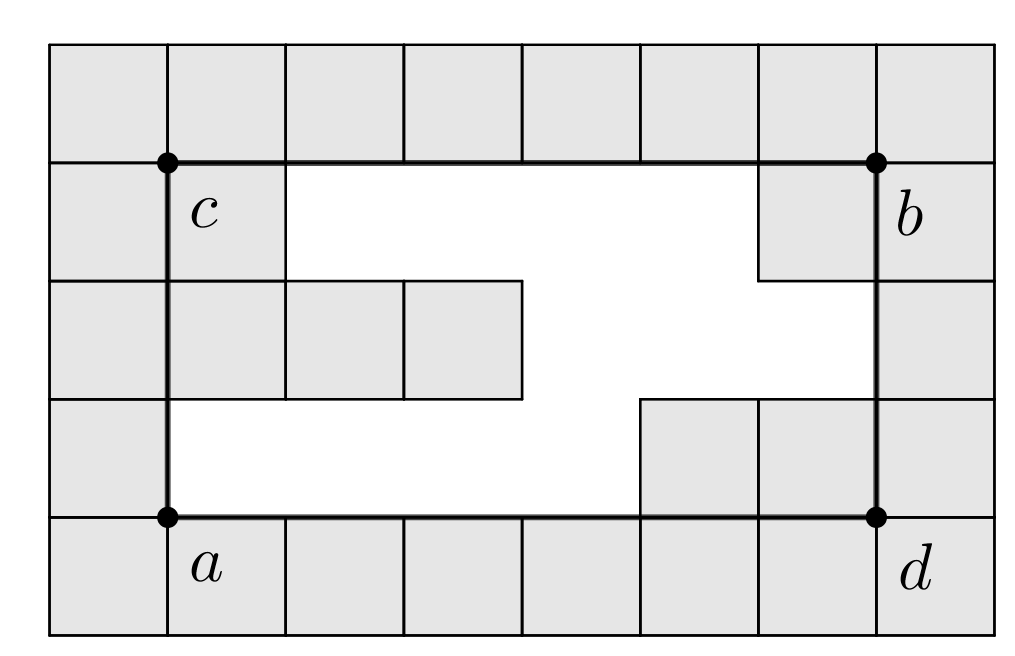}}
    \caption{Polyominoes $\Pc^{c}$.}
    \label{fig:Complement}
\end{figure}

Let $I=[a,b]\subset \NN^{2}$ be a proper interval and let $\Pc_{I}$ denote the
rectangular polyomino determined by $I$.  
If $\Pc\subset \Pc_I$ is a subpolyomino, following~\cite{HibiQureshi2015} we
define the \emph{complement polyomino}
\[
\Pc^{c} \;=\; \Pc_{I} \setminus \Pc .
\]
Whenever $\Pc^{c}$ is a polyomino (for example, when $\Pc$ is convex and does
not meet the boundary of $\Pc_I$), it is non-simple and has exactly one hole
corresponding to the removed region~$\Pc$, see Figure~\ref{fig:Complement} (A).

\begin{Theorem}[\cite{HibiQureshi2015}, Theorem~4.1]\label{Thm:HQ-no-graph}
Let $\Pc\subset\Pc_{I}$ be a simple polyomino and let $\Pc^{c}=\Pc_{I}\setminus
\Pc$.  
Then $I_{\Pc^{c}}$ cannot arise as the toric ideal of the edge ring of any
finite simple graph.
\end{Theorem}

The proof exploits the fact that toric ideals of graphs necessarily contain
binomial $x_{a}x_{b}-x_{c}x_{d}$ (See Figure~\ref{fig:Complement} (B)), while this binomial does not appear in $I_{\Pc^{c}}$.

The main result of~\cite{HibiQureshi2015} shows that, although the graphical
toric interpretation breaks down, the non-simple polyominoes $\Pc^{c}$ obtained
from convex $\Pc$ still have \emph{prime} ideals.

\begin{Theorem}[\cite{HibiQureshi2015}, Theorem~3.1]\label{Thm:HQ-prime}
Let $\Pc\subset \Pc_{I}$ be a convex polyomino such that $\Pc^{c}=\Pc_{I}\setminus \Pc$ is a polyomino.
Then $I_{\Pc^{c}}$ is a prime ideal.
\end{Theorem}

\begin{proof}[Sketch of proof]
Choose the interval $I=[a,b]$ so that $\Pc$ is a convex subpolyomino of
$\Pc_I$, and let $c,d$ be the anti-diagonal corners of $I$ with $b$ and $c$ in
horizontal position.  Using a Gröbner basis criterion for polyomino ideals given in \cite{Qureshi2012}, one first observes
that $x_c$ does not divide the initial term of any binomial in the reduced
Gröbner basis of $I_{\Pc_c}$ with respect to a suitable lex order, \cite[Corollary 2.2]{HibiQureshi2015}.  
Hence $x_c$ is a non-zero divisor on $S_{\Pc_c}/I_{\Pc_c}$.

Localizing at $x_c$ yields an injective map
\[
S_{\Pc^c}/I_{\Pc^c}
\longrightarrow
(S_{\Pc^c}/I_{\Pc^c})_{x_c}
\cong
(S_{\Pc^c})_{x_c}/(I_{\Pc^c})_{x_c}.
\]
The key step is to identify the localized ideal $(I_{\Pc_c})_{x_c}$ with the
polyomino ideal of a \emph{simple} subpolyomino $\Pc'$ of $\Pc^c$, i.e.
\[
(I_{\Pc^c})_{x_c} \;=\; I_{\Pc'} (S_{\Pc^c})_{x_c}.
\]
Since $\Pc'$ is simple, $I_{\Pc'}$ is prime by the results discussed in
Section~\ref{Section:SimplePolyominoes}. Therefore
$(S_{\Pc^c})_{x_c}/(I_{\Pc^c})_{x_c}$ is an integral domain.  
Injectivity of the localization map then implies that
$S_{\Pc^c}/I_{\Pc^c}$ is a domain.
\end{proof}

By \cite[Corollary 2.2]{HibiQureshi2015}, the ideal $I_{\Pc^{c}}$ admits a reduced
quadratic Gröbner basis.  
This yields the following structural properties.

\begin{Theorem}\label{thm:HibiQureshi}
Let $\Pc$ be a convex polyomino such that $\Pc\subset\Pc_{I}$ and
$\Pc^{c}=\Pc_{I}\setminus \Pc$ is a non-simple polyomino.  
Then the coordinate ring $K[\Pc^{c}]$ is a Koszul, normal, Cohen--Macaulay
domain.
\end{Theorem}

\begin{proof}
By \cite[Corollary~2.2]{HibiQureshi2015}, $I_{\Pc^{c}}$ has a squarefree quadratic
Gröbner basis; hence $K[\Pc^{c}]$ is Koszul.  
Normality follows from \cite[Corollary~4.26]{HHO_Book}.  
Cohen--Macaulayness follows from Hochster’s theorem
\cite[Theorem~6.3.5]{Bruns_Herzog}, since toric rings defined by squarefree
initial ideals are Cohen--Macaulay.
\end{proof}

\begin{Remark}
The work of Hibi and Qureshi provided the first infinite family of
\emph{non-simple} polyominoes whose polyomino ideals are prime.  
Their argument is entirely different from the toric description used in the
simple case, relying instead on localization techniques that crucially depend on
the convexity of the removed region~$\Pc$. The method fails when $\Pc$ is simple.
\end{Remark}

It is well known that a binomial prime ideal is a toric ideal (of a suitable toric ring). Therefore, it is natural to ask what is the toric representation of the coordinate ring of the the class of prime non-simple polyominoes provided in \cite{HibiQureshi2015}. In \cite{Shikama2015}, Shikama constructed 
an explicit toric representation for these coordinate rings as follows. 

Let $I=[a,b]$ and $\Pc^c= \Pc_I\setminus \Pc$, where $\Pc$ is a convex polyomino, as before. Suppose that $\Pc$ does not intersect the boundary cells of $\Pc_I$, that is $\Pc^c$ has exactly one hole determined by $\Pc$. In Shikama's setup, let $\Lambda$ be the collection of intervals of $I$ of two types:
\begin{itemize}
  \item[(i)] a special interval $I_e=[a,e]$, where $e$ is the lowest among all
    leftmost outside corners of $\Pc$;
  \item[(ii)] all maximal horizontal and vertical edge intervals of $\Pc$.
\end{itemize}
For each $I\in \Lambda$ introduce a variable $u_I$, and define
\[
\alpha(v)
= \prod_{\substack{I\in\Lambda \\ v\in I}} u_I,
\]
for all $v\in V(\Pc).$ The toric ring associated to $\Pc$ is
\[
T_{\Pc}
= K[\alpha(v) \mid v\in V(\Pc)]
\subset
K[u_I \mid I\in\Lambda].
\]

Let $\varphi \colon S_{\Pc}\to T_{\Pc}$ be the $K$-algebra homomorphism defined by $\varphi(x_v)=\alpha(v)$, and denote the kernel of $\varphi$ by $J_{\Pc}$.

\begin{Theorem}\cite[Theorem 2.3]{Shikama2015}\label{Thm:Shikama-toric}
Let $\Pc^{c}$ be a non-simple polyomino obtained by removing a convex polyomino
$\Pc$ from the rectangle $\Pc_{I}$.  
Then $I_{\Pc}=J_{\Pc}$.  
In particular, $K[\Pc]\cong T_{\Pc}$ is a toric ring and $I_{\Pc}$ is prime.
\end{Theorem}

Shikama proves that $J_{\Pc}$ is generated by quadratic binomials, each
corresponding to an inner $2$-minor of $\Pc$, so that $I_{\Pc}=J_{\Pc}$ holds.

This toric parametrization forms the basis of later developments, including the
zig--zag walk criterion and the characterization of further families of
non-simple polyominoes with prime polyomino ideals.
%%%%%%%%%%%%%%%%%%%%%%%%%%%%%%%%%%%%%%%%%%%%%%%%%%%%%%

\subsection{Zig--zag walks and a necessary condition for primality}
\label{subsec:zigzag}

The next step in the study of non-simple polyominoes is due to Mascia, Rinaldo and
Romeo~\cite{MasciaRinaldoRomeo2020}.  They consider arbitrary non-simple polyominoes and introduce a combinatorial configuration of inner intervals, called a \emph{zig--zag walk}, whose presence always forces the polyomino ideal to be non-prime. 

\begin{Definition}[Zig--zag walk]\label{def:zigzag}
Let $\Pc$ be a polyomino.  A \emph{zig--zag walk} of $\Pc$ is a finite sequence
of distinct inner intervals
\[
\mathcal{W} :\ I_1,\dots,I_{\ell}
\]
with $\ell\ge 2$ such that, for each $i=1,\dots,\ell$, the interval $I_i$ has
either the diagonal corners $v_i,z_i$ and anti-diagonal corners $u_i, v_{i+1}$ or the anti-diagonal corners $v_i,z_i$ and diagonal corners $u_i, v_{i+1}$, and the following
conditions are satisfied:
\begin{enumerate}
  \item[(i)] $I_i\cap I_{i+1}=\{v_{i+1}\}$ for $i=1,\dots,\ell-1$ and
  $I_1\cap I_{\ell}=\{v_1=v_{\ell+1}\}$;
  \item[(ii)] $v_i$ and $v_{i+1}$ lie on the same (horizontal or vertical) edge
  interval of $\Pc$ for all $i$;
  \item[(iii)] for any distinct $i,j$ there is no inner interval $J$ of $\Pc$
  containing both $z_i$ and $z_j$.
\end{enumerate}
Figure \ref{fig:zigzag} shows an example of a zig-zag walk.
\begin{figure}[h]
    \centering
    \includegraphics[scale=0.6]{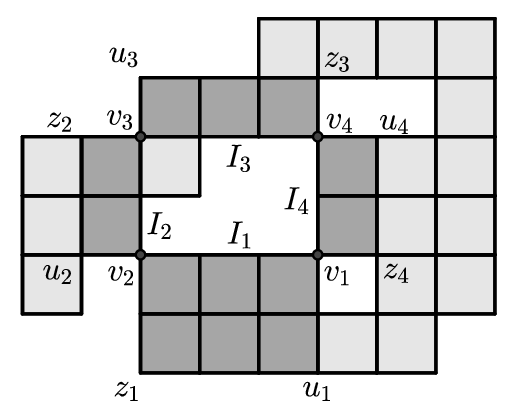}
    \caption{An example of a zig-zag walk.}
    \label{fig:zigzag}
\end{figure}

\end{Definition} 

Condition (ii) forces the sequence of intervals to alternate across edge intervals, which explains the name ``zig--zag''.  From (i) and a simple parity
argument one deduces that the number $\ell$ of intervals in a zig--zag walk is
always even. Moreover, if $v_i$ is the diagonal corner of $I_i$, then $v_{i+1}$ is an anti-diagonal corner of $I_{i+1}$, see \cite[Remark 3.3]{MRR2020}

Each zig--zag walk gives rise to a distinguished binomial. More precisely, if
$\mathcal{W}: I_1,\dots,I_{\ell}$ is a zig--zag walk as above, with notation as
in Definition~\ref{def:zigzag}, set
\begin{equation}\label{eq:zigzag}
f_{\mathcal{W}}
\;=\;
\prod_{k=1}^{\ell} x_{z_k} \;-\;
\prod_{k=1}^{\ell} x_{u_k}.
\end{equation}
Mascia, Rinaldo and Romeo show that $f_{\mathcal{W}}\in J_{\Pc}$ and that it
becomes a zero-divisor modulo $I_{\Pc}$. 

\begin{Proposition}\cite[Proposition 3.5]{MRR2020}\label{prop:zigzag-zerodiv}
Let $\Pc$ be a polyomino and $I_{\Pc}$ its polyomino ideal.  If $\Pc$ admits a
zig--zag walk $\mathcal{W}$, then the elements
\[
x_{v_1},\dots,x_{v_{\ell}},\qquad f_{\mathcal{W}}
\]
are zero divisors in $K[\Pc]$, and
$x_{v_i}f_{\mathcal{W}}\in I_{\Pc}$ for all $i=1,\dots,\ell$.
\end{Proposition}

As a consequence one obtains an intrinsic obstruction to primality.

\begin{Corollary}\cite[Corollary 3.6]{MRR2020}\label{cor:zigzag-nonprime}
Let $\Pc$ be a polyomino.  If $\Pc$ contains a zig--zag walk, then $I_{\Pc}$ is
not a prime ideal.
\end{Corollary}

The authors systematically enumerated non-simple polyominoes up to rank
$14$ and tested primality of their polyomino ideals using \texttt{Macaulay2}.  In this range, the zig--zag obstruction is the only obstruction to primality.

\begin{Theorem}\cite[Proposition 3.9]{MRR2020}\label{thm:rank14}
Let $\Pc$ be a polyomino with $\operatorname{rank}(\Pc)\le 14$.  Then the
following conditions are equivalent:
\begin{enumerate}
  \item $I_{\Pc}$ is a prime ideal;
  \item $\Pc$ contains no zig--zag walk.
\end{enumerate}
\end{Theorem}

The above computations lead to the following general conjecture.

\begin{Conjecture}\cite[Conjecture 4.6]{MRR2020}\label{conj:zigzag}
Let $\Pc$ be a polyomino.  The following conditions are equivalent:
\begin{enumerate}
  \item $I_{\Pc}$ is a prime ideal;
  \item $\Pc$ contains no zig--zag walk.
\end{enumerate}
\end{Conjecture}

%%%%%%%%%%%%%%%%%%%%%%%%%%%%%%%%%%%%%%%%%%%%%%%%%%%%%%%%%%%%%%%%%%%%%%%%%%%%%%%%%%%%%%%%%%%%%%%%%%%%%%%%%%%%%%%%%%%%%%%%%%%%%%%%%%%%%%%%%%%%%%%%%%%%%%%%%%%%%%%%%%%%%%%%%%%%%%%%%%%%%%%%%%%%%%%%%%%%%%%%%%%%%%%%%%%%%%%%%%%%%%%%%%%%%%%%%%%%%%%%%%%%%%%%%%%%%%%%%%%%%%%%%%%%%%%%%%%%%%%%%%%%%%%%%%%%%%%%%%%%%%%%%%%%%%%%%%%%%%%%%%%%%%%%%%%%%%%%%%%%%%%%%%%%%%%%%%%%%%%%%%%%%%%%%%%%%%%%%%%%%%%%%%%%%%%%%%%%%%%%

There are two prominent classes of polyominoes for which the above conjecture holds, namely, closed path (also known as thin cycles) and grid polyominoes.

\subsection*{Grid polyominoes}
Grid polyominoes are non-simple polyominoes obtained by removing a rectangular
grid of pairwise disjoint rectangles from a large rectangle, under strong
alignment conditions on the positions of the holes.  
The following definition is equivalent to \cite[Definition~4.1]{MRR2020}.

\begin{Definition}[Grid polyomino]\label{def:grid}
Let $I=[(1,1),(m,n)]\subset\NN^2$ and let
\[
\Pc \;=\; \Pc_I \setminus \bigcup_{i\in [r],\,j\in[s]} H_{ij},
\]
where each $H_{ij}$ is a rectangle $[a_{ij},b_{ij}]$ strictly contained in $I$,
and the family $\{H_{ij}\}$ satisfies the following alignment conditions:
\begin{enumerate}
  \item for fixed $i$, all $H_{ij}$ have the same $x$-coordinates of their
  vertical sides;
  \item for fixed $j$, all $H_{ij}$ have the same $y$-coordinates of their
  horizontal sides;
  \item consecutive holes in a row or column are separated by exactly one layer
  of cells.
\end{enumerate}
Then $\Pc$ is called a grid polyomino. See Figure \ref{fig:grid} for an example of a grid polyomino.

\begin{figure}[h]
    \centering
    \includegraphics[scale=0.6]{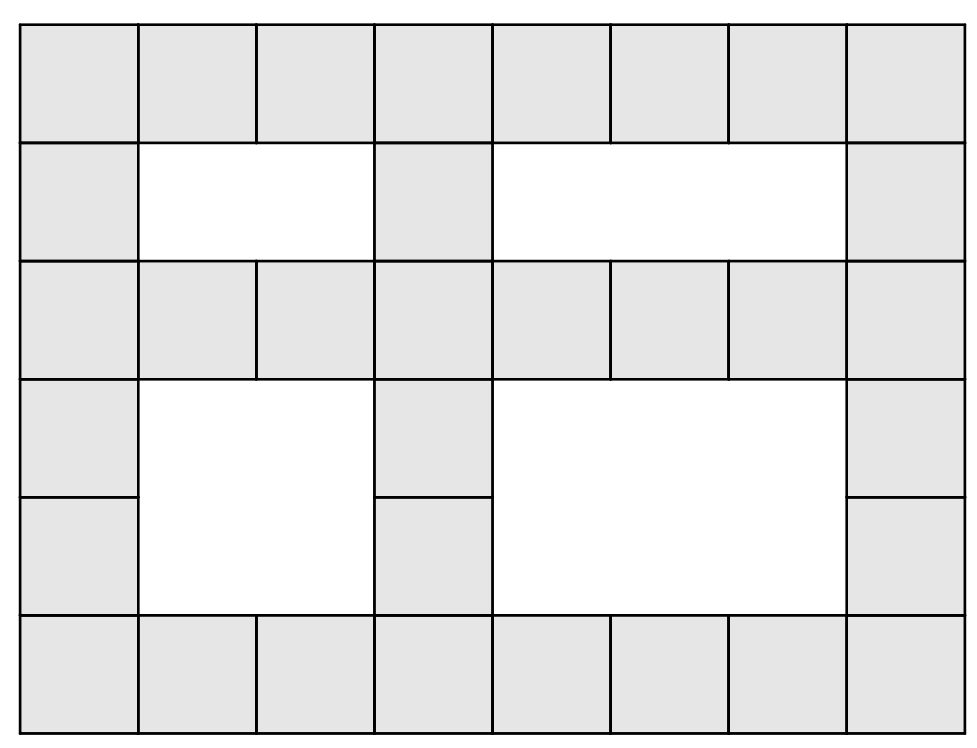}
    \caption{A grid polyomino.}
    \label{fig:grid}
\end{figure}
\end{Definition} 

Grid polyominoes are non-simple and may contain many holes, but their
highly regular structure prevents the existence of zig--zag walks. Figure~\ref{fig:grid} displays an example of a grid polyomino. To show that grid polyominoes are prime, the authors in \cite{MRR2020} generalized Shikama’s construction of toric rings \cite{Shikama2015} to an arbitrary non-simple polyomino.

Let $\Pc$ be a polyomino.  
Denote by $\{V_i\}_{i\in I}$ and $\{H_j\}_{j\in J}$ the sets of maximal
vertical and horizontal edge intervals of $\Pc$, respectively, and by
$H_1,\dots,H_r$ the holes of $\Pc$.  
For each hole $H_k$ let $e_k=(i_k,j_k)$ be its lower-left corner, and set
\[
F_k \;=\; \{(i,j)\in V(\Pc) : i\le i_k,\ j\le j_k\}.
\]
Introduce variables $v_i$ for $V_i$, $h_j$ for $H_j$ and $w_k$ for $F_k$, and define
\[
\psi \colon V(\Pc)\longrightarrow
T'=K[v_i,h_j,w_k \mid i\in I,\ j\in J,\ k=1,\dots,r]
\]
by
\[
\psi(a)\;=\;
\prod_{a\in H_i\cap V_j} h_i v_j \cdot
\prod_{a\in F_k} w_k.
\]
Let
\[
T'_{\Pc} \;=\; K[\psi(a)\mid a\in V(\Pc)]\subset T'
\]
be the associated toric ring, and let $J'_{\Pc}$ denote the kernel of the
surjective homomorphism
\[
\Phi' \colon S_{\Pc}=K[x_a : a\in V(\Pc)] \longrightarrow T'_\Pc,\qquad
\Phi'(x_a)=\psi(a).
\]

By construction $J'_{\Pc}$ is a prime binomial ideal containing the polyomino
ideal $I_{\Pc}$.  
A key observation in \cite[Lemma~3.1]{MRR2020} is that the quadratic part of
$J'_{\Pc}$ coincides with $I_{\Pc}$.  
Thus $J'_{\Pc}$ may be viewed as a canonical toric enlargement of $I_{\Pc}$.

It remains an open and challenging problem to describe the elements of
$J'_{\Pc}\setminus I_{\Pc}$.  
It is observed in \cite[Remark~3.7]{MRR2020} that if $\Pc$ contains a zig--zag
walk, then its associated binomial defined in~\eqref{eq:zigzag} belongs to
$J'_\Pc$.  
However, not every element of $J'_{\Pc}\setminus I_{\Pc}$ is of this form, as
illustrated in \cite[Example~3.8]{MRR2020}.

Given a grid polyomino $\Pc$, using a redundancy criterion for binomials of degree at least $3$, generalizing a lemma used by Shikama in \cite{Shikama2015}, the authors in \cite{MRR2020} proved that every irredundant binomial in $J_{\Pc}$ has degree $2$. In particular, all higher-degree binomials in
$J'_{\Pc}$ are redundant.

\begin{Theorem}\cite[Theorem]{MRR2020}\label{thm:grid-prime}
Let $\Pc$ be a grid polyomino.  Then $I_{\Pc}=J_{\Pc}$ and, in particular,
$I_{\Pc}$ is a prime ideal.
\end{Theorem}
%%%%%%%%%%%%%%%%%%%%%%%%%%%%%%%%%%%%%%%%%%%%%%%%%%%%%%%%%%%%%%%%%%%%%%%%%%%%%%%%%%%%%%%%%%%%%%%%%%%%%%%%%%%%%%%%%%%%%%%%%%%%%%%%%%%%%%%%%%%%%%%%%%%%%%%%%%%%%%%%%%%%%%%%%%%%%%%%%%%%%%%%%%%%%%%%%%%%%%%%%%%%%%%%%%%%%%%%%%%%%%%%%%%%%%%%%%%%%%%%%%%%%%%%%%%%%%%%%%%%%%%%%%%%%%%%%%%%%%%%%%%%%%%%%%%%%%%%%%%%%%%%%%%%%%%%%%%%%%%%%%%%%%%%%%%%%%%%%%%%%%%%%%%%%%%%%%%%%%%%%%%%%%%%%%%%%%%%%%%%%%%%%%%%%%%%%%%%%%%%%%%%%%%%%%%%%%%%%%%%%%%%%%%%%%%%%%%%%%%%%%%%%%%%%%

\subsection*{Closed path polyominoes.} 
Closed path polyominoes were introduced and studied in detail by Cisto and Navarra~\cite{CistoNavarra2021}.  
Informally, a closed path is a non-simple polyomino with exactly one hole, obtained by arranging the cells in a cyclic fashion so that they form a loop.  
The following definition formalizes this notion.
\begin{Definition}
Let $A_1, A_2, \ldots, A_n$ be a sequence of distinct cells with $n>5$.  
The sequence is called a \emph{closed path} if it is a polyomino and satisfies the following conditions:
\begin{enumerate}
\item consecutive cells share an edge; that is, $A_i \cap A_{i+1}$ is a common edge of both cells for all $i=1,\ldots,n$, and 
\item if $i\in\{1,\ldots,n\}$ and $j\notin\{i-2,\, i-1,\, i,\, i+1,\, i+2\}$, then $A_i \cap A_j = \emptyset$.
\end{enumerate}
Here boundary indices are interpreted cyclically by setting $A_{-1}=A_{n-1}$, $A_0=A_n$, $A_{n+1}=A_1$, and $A_{n+2}=A_2$. See Figure \ref{fig:closed path} for an example of a closed path polyomino, on the left, and a non-closed path polyomino, on the right.

\begin{figure}[h]
    \centering
    \subfloat{\includegraphics[width=0.3\linewidth]{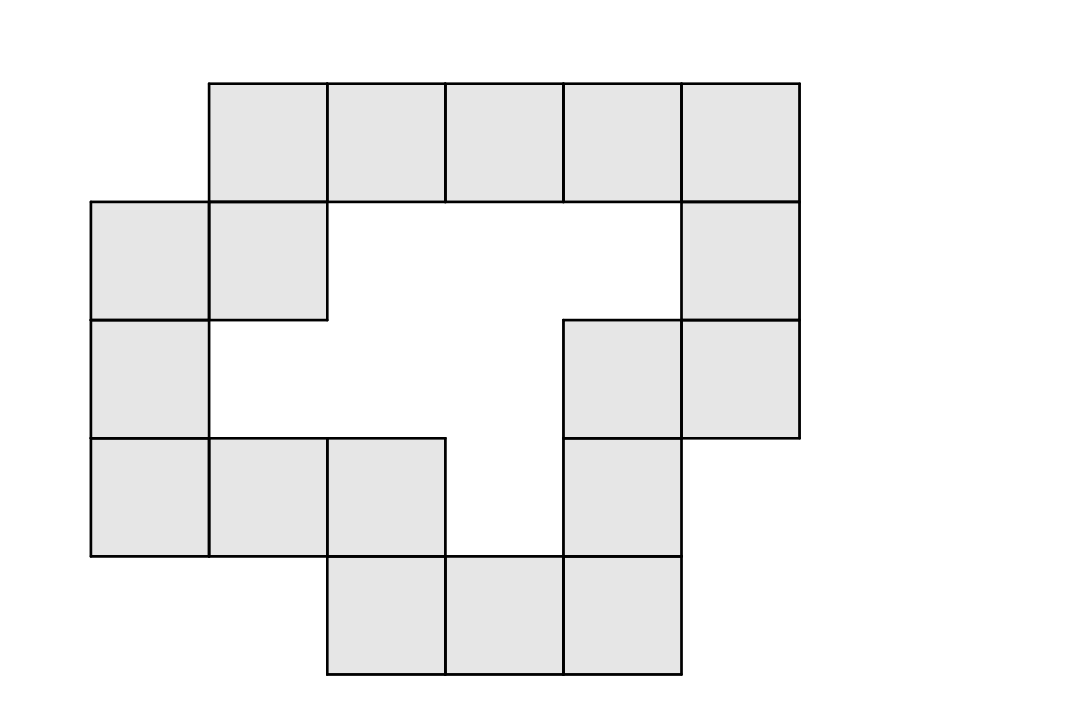}}\qquad
    \subfloat{\includegraphics[width=0.3\linewidth]{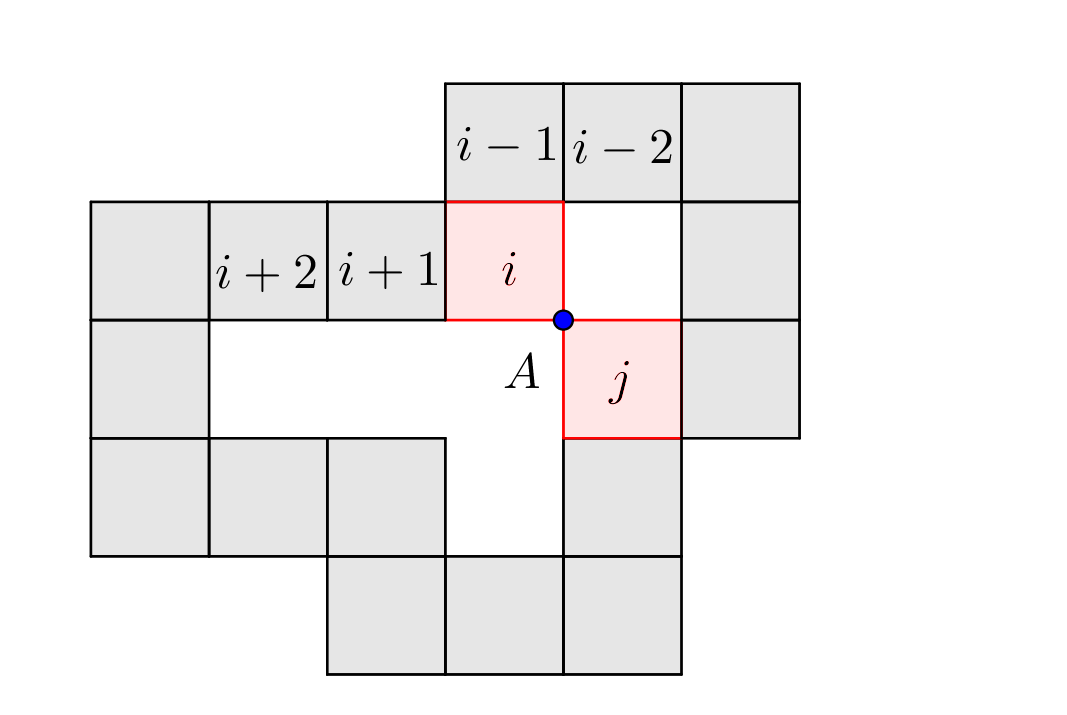}}
    \caption{A closed path and a non-closed path.}
    \label{fig:closed path}
\end{figure}
\end{Definition}

It follows from the definition that the closed paths are thin polyominoes and can be regarded as ``necklaces of cells'' surrounding a hole. To give a complete characterization of prime closed path polyominoes, the following possible configurations within them play an important role.

\begin{figure}[h]
    \centering
    \subfloat[]{\includegraphics[width=0.3\linewidth]{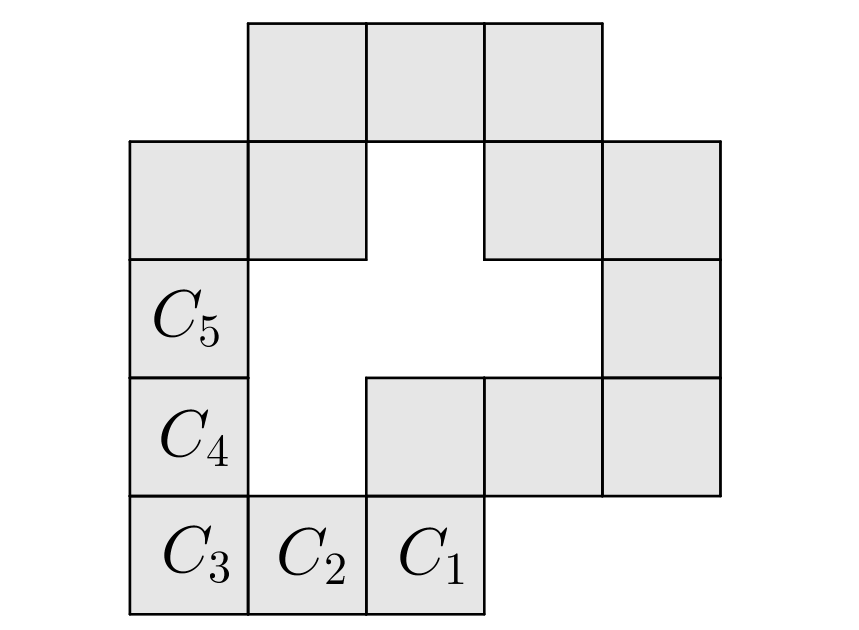}}\qquad
    \subfloat[]{\includegraphics[width=0.3\linewidth]{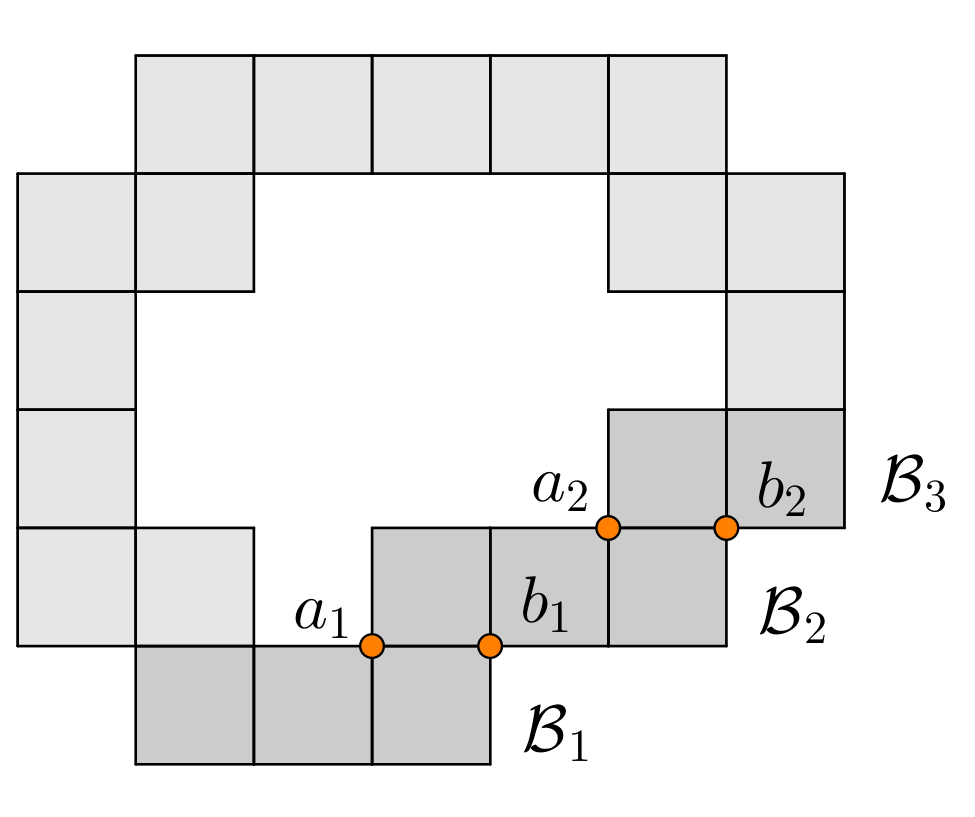}}
    \caption{A closed path with an $L$-configuration and another one with a ladder of three steps.}
    \label{fig:L conf and ladder}
\end{figure}

\begin{itemize}
  \item An \emph{$L$-configuration}, consisting of a path of five cells
        $C_1,\dots,C_5$ such that $C_1,C_2,C_3$ and $C_3,C_4,C_5$ form two
        orthogonal blocks; see Figure~\ref{fig:L conf and ladder} (A).
\item A \emph{ladder with at least three steps}, namely a sequence of
maximal horizontal or vertical blocks arranged alternately so that each block meets the next in a single vertex, and the orientation switches at every step; see Figure~\ref{fig:L conf and ladder} (B).
\end{itemize}

In \cite{CistoNavarra2021}, a complete characterization of prime closed path polyominoes, and an affirmative answer to Conjecture\ref{conj:zigzag} is provided. Below we briefly sketch these ideas. 

\begin{Theorem}\label{thm:closedpath}
Let $\Pc$ be a closed path polyomino.  
Then the following conditions are equivalent:
\begin{enumerate}
  \item $I_{\Pc}$ is prime;
  \item $\Pc$ contains no zig--zag walks;
  \item $\Pc$ contains either an $L$-configuration or a ladder with at least three steps.
\end{enumerate}
\end{Theorem}

\begin{proof}[Sketch of proof]
The implication (1) $\Rightarrow$ (2) follows from 
Corollary~\ref{cor:zigzag-nonprime}.  

For (2) $\Rightarrow$ (3), assume that $\Pc$ is a closed path polyomino that contains no zig–-zag walks.  
A detailed combinatorial analysis shows that this geometric restriction forces the presence of a controlled local configuration around the hole: namely, $\Pc$ must contain either an $L$–configuration or a ladder with at least three steps.  
This is established in \cite[Proposition~6.1]{CistoNavarra2021}.

To prove (3) $\Rightarrow$ (1), one constructs an explicit toric parametrization of $K[\Pc]$ using a refinement of Shikama’s toric method for non-simple polyominoes \cite{Shikama2015}, introducing a single additional ``hole variable.’’  
Assume that $\Pc$ contains either an $L$–configuration or a ladder with at least three steps. Let $A$ be the set of vertices singled out by this configuration (the vertices
inside the $L$–shape in the first case, or the vertices
inside the ladder in the
second). As usual, let $\{V_i\}_{i\in I}$ and $\{H_j\}_{j\in J}$ be the sets of maximal vertical and horizontal edge intervals of $\Pc$, with associated variables $\{v_i\}_{i\in I}$ and $\{h_j\}_{j\in J}$, and let $w$ be an extra variable.  

Define a map
\[
\alpha : V(\Pc)\longrightarrow K[\{v_i,h_j,w\}_{i\in I,\, j\in J}],
\]
by
\[
\alpha(r)= v_i h_j\, w^k,
\]
whenever $V_i\cap H_j=\{r\}$, where $k=0$ if $r\notin A$ and $k=1$ if $r\in A$.  

Let $T_{\Pc}=K[\alpha(v) : v\in V(\Pc)]$ and let $J$ be the kernel of the surjective $K$–algebra homomorphism
\[
\phi : S_{\Pc} \longrightarrow T_{\Pc}, \qquad \phi(x_v)=\alpha(v),
\]
for all $v\in V(\Pc)$.  
It is proved in \cite[Theorem~4.2 and Theorem~5.2]{CistoNavarra2021} that $I_{\Pc}=J$ by expressing binomials in $J$ as combinations of quadratic binomials corresponding exactly to the inner $2$–minors of $\Pc$.  
Hence $I_{\Pc}$ is the defining toric ideal of $T_{\Pc}$, and in particular it is prime.
\end{proof}

\begin{figure}[h]
    \centering
    \includegraphics[width=0.3\linewidth]{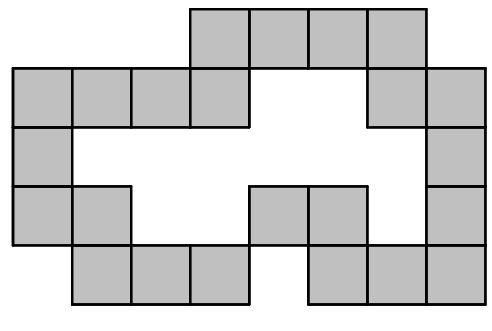}
    \caption{A weakly closed path.}
    \label{fig:weakly closed}
\end{figure}

In subsequent work with Utano \cite{CNU2}, this strategy is
extended to \emph{weakly} closed path polyominoes (see Figure \ref{fig:weakly closed} for an example), giving further evidence that
zig--zag walks play a central role in the characterization of prime polyomino ideals.

%%%%%%%%%%%%%%%%%%%%%%%%%%%%%%%%%%%%%%%%%%%%%%%%%%%%%%%%%%%%%%%%%%%%%%%%%%%%%%%%%%%%%%%%%%%%%%%%%%%%%%%%%%%%%%%%%%%%%%%%%%%%%%%%%%%%%%%%%%%%%%%%%%%%%%%%%%%%%%%%%%%%%%%%%%%%%%%%%%%%%%%%%%%%%%%%%%%%%%%%%%%%%%%%%%%%%%%%%%%%%%%%%%%%%%%%%%%%%%%%%%%%%%%%%%%%%%%

\subsection{Quadratic Gr\"obner bases and primality}
Mascia, Rinaldo and Romeo, in~\cite{MasciaRinaldoRomeo2020}, developed an 
interesting criterion to establish the primality of polyomino ideals through 
the existence of \emph{quadratic Gr\"obner bases} of $I_{\Pc}$.  
Their method is purely combinatorial and is based on Gröbner bases with 
respect to certain graded reverse lexicographic orders naturally induced by 
the geometry of $\Pc$. They introduced eight graded reverse lexicographic orders on \(S_{\Pc}\), each obtained by reading the vertices of \(\Pc\) according to one of the eight directions in the plane, and for each such order \(<_{i}\), they give an explicit combinatorial criterion ensuring that the inner 2-minors form a reduced Gröbner basis of \(I_{\Pc}\). These criteria are expressed through eight local geometric configurations \((\pi_{1}),\dots,(\pi_{8})\) (see~\cite[Definition 3.3 and Table 2]{MasciaRinaldoRomeo2020} ).

Given any vertex \(v\in V(\Pc)\), they refine a chosen order \(<_{i}\) to a new order \(<_{i,v}\) by declaring \(x_{v}\) to be the smallest variable while keeping the relative order of all other variables. Their central structural result shows:
\begin{itemize}
    \item If the inner 2-minors form a reduced Gröbner basis of \(I_{\Pc}\) for some \(<_{i}\),
    \item and if every vertex \(v\) fails at least one of the conditions \(\pi_{k}\) corresponding to that order,
\end{itemize}
then the inner 2-minors also form a reduced Gröbner basis with respect to the modified order $<_{i,v}$ for every vertex $v$.
Consequently, $I_{\Pc} = (I_{\Pc} : u)$ for every monomial $u \in S_\Pc$.
Since $(I_{\Pc} : u)$ equals the lattice ideal of the saturated lattice $\Lambda$ whose basis corresponds to the cells in $\Pc$, it follows that $I_{\Pc}$ itself is prime.

As an application of this approach, it is shown that $I_{\Pc}$ is prime when $\Pc$ correspond to {\em subgrid polyominoes}, obtained by removing certain cells from a grid polyomino while preserving connectivity.

Moreover, by using \cite[Corollary 3.3]{MasciaRinaldoRomeo2020}, Koley, Kotal and Veer obtained the following in \cite{KoleyKotalVeer2024}.

\begin{Proposition}\cite[Proposition 5.4.]{KoleyKotalVeer2024} 
 Let $\cP$ be a thin polyomino such that if two maximal inner intervals of $\cP$ intersects, then there intersection is a cell. Then $I_{\cP}$ is a prime ideal.
\end{Proposition} 

% This provides strong partial evidence for the general conjecture that the absence of zig-–zag walks should characterize primality. 

%%%%%%%%%%%%%%%%%%%%%%%%%%%%%%%%%%%%%%%%%%%%%%%%%%%%%%%%%%%%%%
%%%%%%%%%%%%%%%%%%%%%%%%%%%%%%%%%%%%%%%%%%%%%%%%%%%%%%%%%%%%%%%%%%%%%%%%%%%%%%%%%%%
\section{Radicality and Primary Decomposition of Polyomino Ideals}
\label{Section:Radicality}
Although the theory of polyomino ideals has developed extensively in recent years, 
most work has focused on understanding when such ideals are prime.  
In contrast, much less is known about radicality and about the primary 
decompositions of polyomino ideals that are not prime.  
At present, no example is known of a polyomino ideal that is prime but not radical, 
and it is natural to ask:

\begin{Question}
Are all polyomino ideals radical?
\end{Question}

A standard criterion for radicality asserts that if an ideal has a squarefree
initial ideal with respect to some monomial order, then the ideal itself is
radical. The first systematic study of radicality for a non-prime class of
polyominoes using this approach was carried out in \cite{YonatanHamonangan2019}.  
They introduced a class of non-simple polyominoes, called \emph{cross polyominoes}, defined as a union of two
rectangles satisfying certain intersection conditions
(see~\cite[Definition~3.2]{YonatanHamonangan2019}).  
Typically, cross polyominoes are non-prime. The authors in \cite{YonatanHamonangan2019} proved that if
the intersection consists of a single cell, then the associated polyomino ideal
is radical.  
Their proof constructs a Gröbner basis whose initial ideal is squarefree with
respect to a suitable monomial order, and then applies the general radicality
criterion.

\begin{figure}[h]
    \centering
    \includegraphics[scale=1]{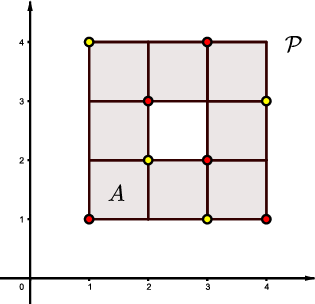}
    \caption{A closed path.}
    \label{fig:Graver}
\end{figure}

As discussed in Section~\ref{Section:SimplePolyominoes}, if $\Pc$ is a simple
polyomino, then $I_{\Pc}$ admits a squarefree initial ideal.  
The same phenomenon occurs for most known classes of non-simple \emph{prime}
polyominoes.  
This naturally raises the question of whether polyomino ideals admit a
\emph{squarefree universal} Gröbner basis.  
A negative answer was given in~\cite[Remark~16]{CistoNavarra2021}, where an
explicit example is provided of a polyomino shown in Figure~\ref{fig:Graver}, whose universal Gröbner basis
contains binomials with non-squarefree monomials, as the binomial $f=x_{11}x_{23}x_{32}x_{34}x_{41}-x_{14}x_{22}x_{31}^2x_{43}$ attached to the vertices in red and yellow.

Such examples show that, even when a polyomino ideal is radical, its universal Gröbner basis need not consist of squarefree binomials.

Very recently, Koley, Kotal, and Veer initiated the study of radicality in connection with the Knutson property in \cite{KoleyKotalVeer2024}. In fact, in arbitrary characteristic, Knutson ideals have squarefree initial ideals and are therefore radical.

\subsection*{Polyocollections and primary decompositions}

A recent contribution to the study of radicals and primary decompositions of
polyomino ideals is due to Cisto, Navarra and Veer~\cite{CistoNavarraVeer2024}. They introduced the notion of \emph{polyocollections} as a combinatorial
generalization of collections of cells.  

\begin{Definition}
    Let $\cC$ be a collection of intervals in $\mathbb{Z}^2$. We say that $\cC$ is a polyocollection if for all $I, J \in \cC$ with $I \ne J$, we have $I \cap J \ne \emptyset$ and one of the following holds:
\begin{enumerate}
    \item $I \cap J$ is a common edge\footnote{If $[a,b]$ is a proper interval, and $c,d$ are its antidiagonal vertices, then $E([a,b])=\{[a,c],\, [a,d],\, [b,d],\, [b,c]\}.$} of $I$ and $J$.
    \item For all $F \in E(I)$ and for all $G \in E(J)$, we have $\lvert F \cap G \rvert \le 1$.
\end{enumerate}
\end{Definition}

\begin{figure}[h!]
		\subfloat[]{\includegraphics[scale=0.6]{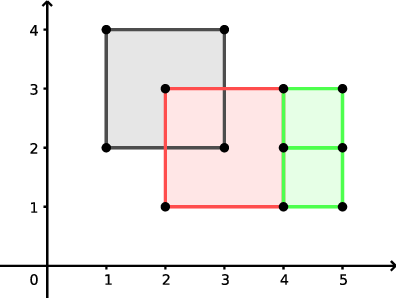}}\qquad \qquad
		\subfloat[]{\includegraphics[scale=0.6]{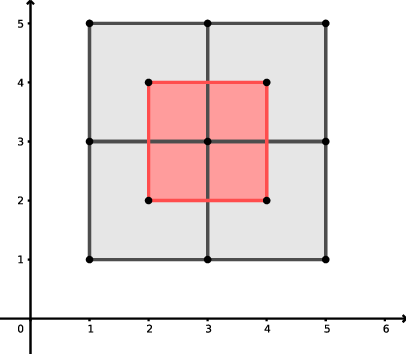}}
		\caption{A polyocollection on left and a non-polyocollection on right.}
		\label{img1}
\end{figure}

For example, the collection
\[
\mathcal{C}_1=\{[(1,1),(3,3)],\,[(1,3),(3,5)],\,[(3,1),(5,3)],\,[(3,3),(5,5)],\,[(2,2),(4,4)]\},
\]
displayed in Figure~\ref{img1}\,(A), is a polyocollection.  
In contrast, 
\[
\mathcal{C}_2=\{[(1,2),(3,4)],\,[(2,1),(4,3)],\,[(4,1),(5,2)],\,[(4,2),(5,3)]\},
\]
displayed in Figure~\ref{img1}\,(B), is not a polyocollection.

To each polyocollection $\cC$, the authors associate a binomial ideal $I_\cC$, extending in a natural way the classical construction of polyomino ideals. They provided a unified framework to primary decomposition, in terms of the so-called \emph{admissible sets} and the lattice ideals of suitable sub-collections derived from $\cC$. In particular, admissible sets are defined as follows:

\begin{Definition}
    A subset $X \subseteq V(\cC)$ is called an \emph{admissible set} of $\cC$ if, for every inner interval $I$ of $\cC$, either $X \cap V(I) = \emptyset$ or $X \cap V(I)$ contains the boundary of an edge of $I$.
\end{Definition}

Using this notion and several algebraic techniques coming from \cite{HHMoradi} and from lattice ideal theory, they prove the following result:

\begin{Theorem}\cite[Proposition 3.12, Theorem 3.13]{CistoNavarraVeer2024}
    Let $\cP$ be a polyomino (or, more generally, a polyocollection). Then, for every minimal prime ideal $\mathfrak{p}$ of $I_\cP$, there exists an admissible set $X$ such that 
    \[
    \mathfrak{p} = J_X := (\{x_a \mid a \in X\}) + L_{\cP^{(X)}},
    \]
    where $L_{\cP^{(X)}}$ is the lattice ideal of the polyocollection $  \cP^{(X)}$ having the set $\{ I \text{ inner interval of } \cC \mid V(I) \cap X = \emptyset \}$ as set of inner intervals.\\    
    Hence,
    \[
    \sqrt{I_\cP} = \bigcap_X J_X,
    \]
    where the intersection runs over all admissible sets $X$ of $\cP$.
\end{Theorem}

For our purposes, the key point is that this framework provides the first systematic method for describing a primary decomposition of non-prime polyomino ideals. In particular, the authors completely determine the minimal primary decomposition of non-prime closed path polyominoes.  

If $\cP$ is a closed path containing a zig–zag walk (equivalently, $I_\cP$ is not prime by Theorem~\ref{thm:closedpath}), then \cite[Theorem 4.19]{CistoNavarraVeer2024} shows that
\[
I_\cP = \mathfrak{p}_1 \cap \mathfrak{p}_2,
\]
where both $\mathfrak{p}_1$ and $\mathfrak{p}_2$ are binomial prime ideals of height $|\cP|$. In particular, $\mathfrak{p}_1$ is the toric ideal appearing in Mascia–Rinaldo–Romeo~\cite{MasciaRinaldoRomeo2020}, while $\mathfrak{p}_2$ is a combinatorially defined monomial–binomial ideal (look at \cite[Notation 4.7]{CistoNavarraVeer2024}).

\begin{figure}[h]
    \centering
    \subfloat[]{\includegraphics[scale=0.8]{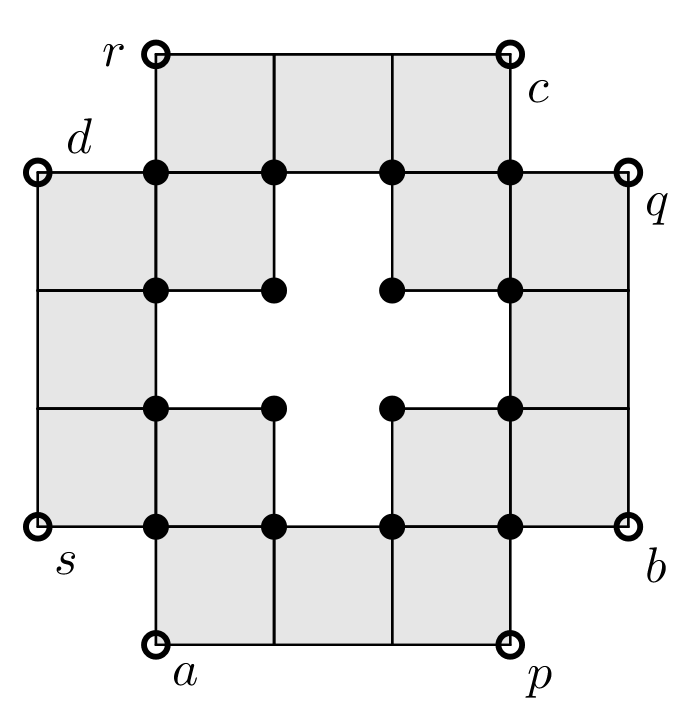}}\qquad\qquad
    \subfloat[]{\includegraphics[scale=0.7]{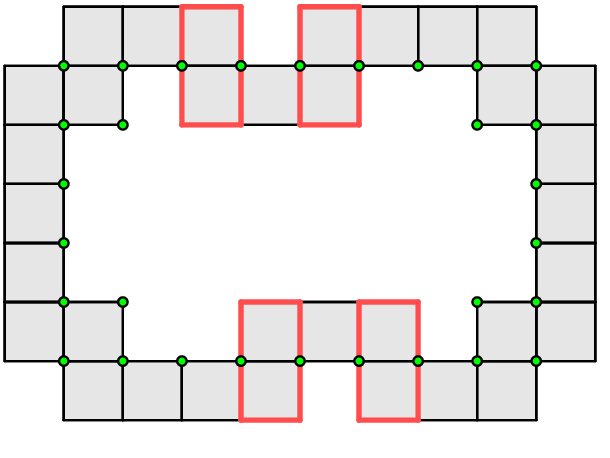}}
    \caption{Non-prime closed paths.}
    \label{fig:prim decomp}
\end{figure}

For instance, the non-prime closed path in Figure~\ref{fig:prim decomp} (A) has
\[
\mathfrak{p}_1 = I_\cP + (x_a x_b x_c x_d - x_p x_q x_r x_s), 
\qquad
\mathfrak{p}_2 = (x_a : a \text{ is a black point}).
\]

The non-prime closed path in Figure~\ref{fig:prim decomp} (B) has
\[
\mathfrak{p}_1 = I_\cP + (f_\cW : \cW \text{ is a zig-zag walk of } \cP),
\]
\[
\mathfrak{p}_2 = (x_a : a \text{ is a green point}) + (\text{binomials attached to the red intervals}).
\]

%Thus every non-prime closed path polyomino is \emph{unmixed}, and its minimal primes admit transparent combinatorial descriptions.

These results provide the strongest evidence so far that the primary
decomposition of polyomino ideals is governed entirely by combinatorics.  
Describing the \textit{minimal} primary decomposition of an arbitrary non-prime polyomino in terms of the combinatorics of the polyomino itself remains a challenging and largely open problem.

\section{Hilbert-Poincar\'e series and rook polynomial theory}
\label{Section:RookPoly}

A recent line of research has discovered a novel connection showing that the Hilbert--Poincar\'e series of $K[\cP]$ is closely related to the rook polynomial and one of its variants, namely the switching rook polynomial. To describe this connection, we first recall the definitions of the $h$-polynomial and the Castelnuovo--Mumford regularity of a graded ideal.\\
Let $R = K[x_1, \dots, x_n]$ be a polynomial ring over a field $K$, and let $I$ be a homogeneous ideal of $R$. The \emph{Castelnuovo--Mumford regularity} (or simply \emph{regularity}) of $I$ is defined by
\[
\reg(I) = \max\{\, j \mid \beta_{i,\, i+j} \neq 0 \text{ for some } i \,\}.
\]
where $\beta_{i,j}$ denotes the \emph{$(i,j)$-th graded Betti number} of $I$. Moreover, one has $\reg(R/I)=\reg(I)-1$. \\
The quotient $R/I$ has a natural grading as a $K$-algebra, that is, $R/I = \bigoplus_{k \in \mathbb{N}} (R/I)_k.$ The associated formal power series $\mathrm{HP}_{R/I}(t) = \sum_{k \in \mathbb{N}} \dim_K (R/I)_k\, t^k$ is called the \emph{Hilbert--Poincar\'e series} of $R/I$.  
By the Hilbert--Serre Theorem, there exists a unique polynomial $h(t) \in \mathbb{Z}[t]$ such that $h(1) \neq 0$ and
\[
\mathrm{HP}_{R/I}(t) = \frac{h(t)}{(1 - t)^d},
\]
where $d$ is the Krull dimension of $R/I$. The polynomial $h(t)$ is called the \emph{$h$-polynomial} of $R/I$. Furthermore, if $R/I$ is Cohen--Macaulay, then $\deg h(t)=\reg(R/I).$

We now recall the notion of non-attacking rooks on a collection of cells.

% The placement of non-attacking rooks on a skew diagram corresponds to the enumeration of permutations with certain restrictions; this idea was introduced by Kaplansky and Riordan~\cite{KR} and further developed by Riordan~\cite{R}.  
% For a comprehensive treatment of permutations with forbidden positions, we refer the reader to Stanley~\cite[Chapter~2]{S1}.

\subsection{Rook polynomial.}	
\begin{figure}[h]
		\centering
		\subfloat[Attacking rooks]{\includegraphics[scale=0.8]{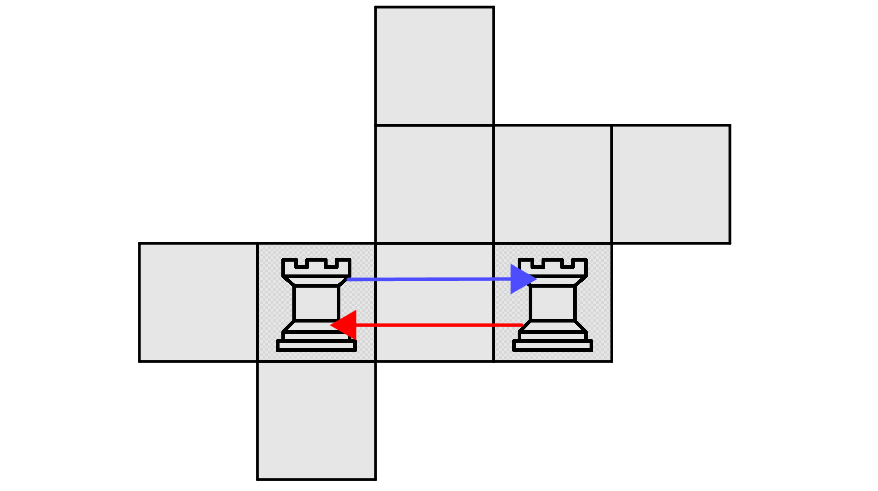}}\quad
		\subfloat[Non-attacking rooks]{\includegraphics[scale=0.8]{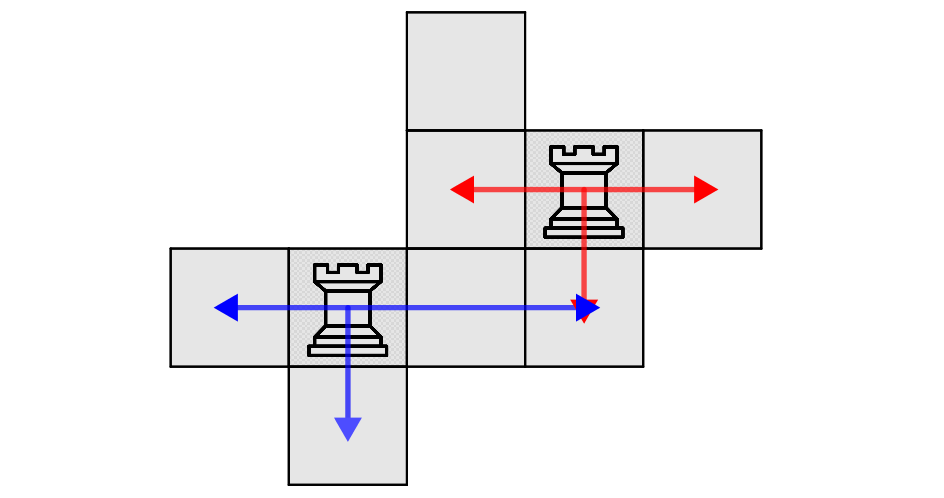}}	
		\caption{Positions of two rooks in a polyomino.}
		\label{Figure: exa attacking rooks}
	\end{figure}

Let $\cP$ be a collection of cells. Two rooks $R_1$ and $R_2$ are in \textit{attacking position} or \textit{attacking rooks} in $\cP$ if there exists two cells $A_1$ and $A_2$ of $\cP$ in horizontal or vertical position such that $R_1$ and $R_2$ are placed in $A_1$ and $A_2$, respectively, and $[A_1,A_2]$ is contained in $\cP$. In contrast, two rooks are in \textit{non-attacking position} or \textit{non-attacking rooks} in $\cP$ if they are not in attacking position. For instance see Figure \ref{Figure: exa attacking rooks}.

\begin{figure}[h]
		\centering
		\includegraphics[scale=0.85]{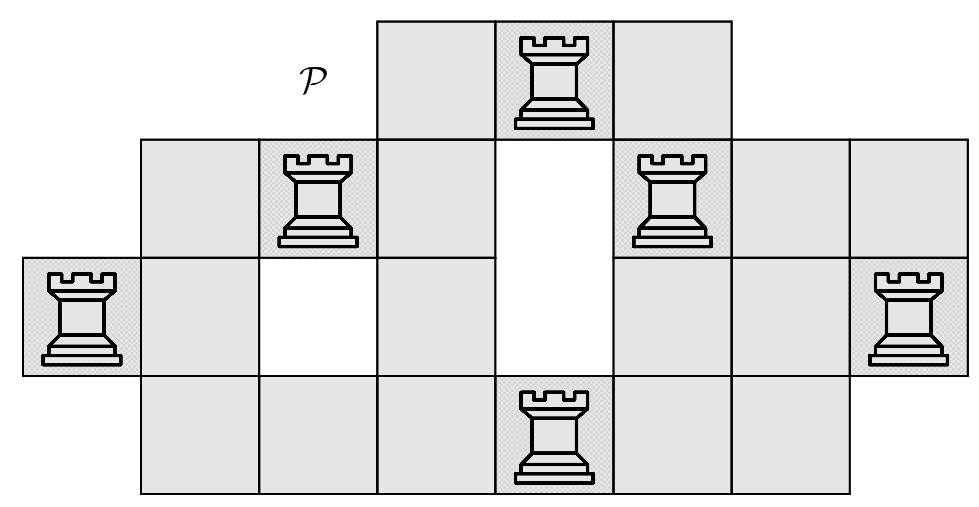}
		\caption{An example of a $6$-rook configuration in $\cP$.}
		\label{Figura:esempio rook configuration}
	\end{figure}

	 A \textit{$j$-rook configuration} in $\cP$ is a set of $j$ rooks arranged in non-attacking positions within $\cP$, where $j\geq 0$; for convention, the $0$-rook configuration is $\emptyset$.  Figure~\ref{Figura:esempio rook configuration} shows a 6-rook configuration. We say that a $j$-rook configuration in $\cP$ is \textit{maximal} if there does not exist any $k$-rook configuration in $\cP$, with $k > j$, that properly contains it.

\begin{figure}[h]
	 	\centering
	 	\includegraphics[width=0.23\textwidth]{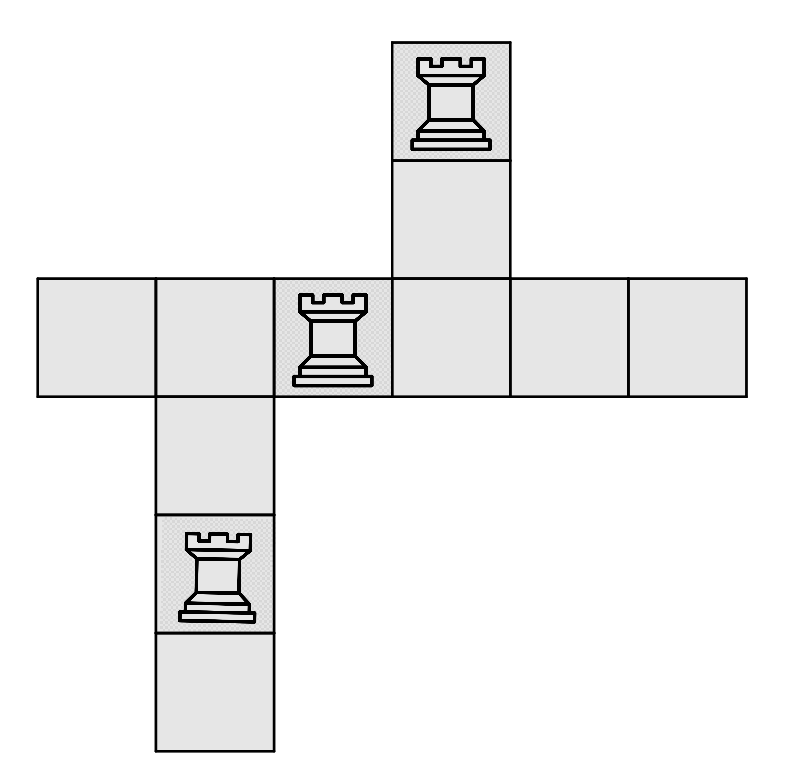}
	 	\caption{Polyomino}
	 	\label{fi:a pol}
	 \end{figure}

	 The \textit{rook number} $r(\cP)$ is the maximum number of rooks that can be placed in $\cP$ in non-attacking positions. We denote by $ \mathcal{R}(\mathcal{P},k)$ the set of all $k$-rook configurations in $\mathcal{P}$ and set $ r_k = \vert \mathcal{R}(\mathcal{P},k) \vert $, for all $ k \in \{0, \dots, r(\mathcal{P})\} $ (with the convention $ r_0 = 1 $). The rook polynomial of $\mathcal{P}$ is the polynomial in $\mathbb{Z}_{> 0}[t]$ defined as $$ r_{\mathcal{P}}(t) = \sum_{k=0}^{r(\mathcal{P})} r_k t^k. $$
	  For instance, the polyomino in Figure \ref{fi:a pol} has $r_{\cP}(t)=1+11t+31t^2+24t^3$ and $r(\cP)=3$.
	 
  In Combinatorial Commutative Algebra, the significance of the rook number and the rook polynomial of a collection of cells $\cP$ arises from their connection with the Castelnuovo–Mumford regularity and the $h$-polynomial of $K[\cP]$, respectively. This connection was first established by Ene et al.\ in 
\cite{EneHerzogQureshiRomeo2021}, where they studied algebraic invariants of $L$-convex polyominoes.

%%%%%%%%%%%%%%%%%%%%%%%%%%%%%%%%%%%%%%%%%%%%%%%%%%%%%%%%%%%%%%%%%%%%%%%%%%%%%%%%%%%%%%%%%%%%%%%%%%%%%%%%%%%%%%%%%%%%%%%%%%%%%%%%%%%%%%%%%%%%%%%%%%%%%%%%%%%%%%%%%%%%%%%%%%%%%%%%%%%%%%%%%%%%%%%%%%%%%%%%%%%%%%%%%%%%%%%%%%%%%%%%%%%%%%%%%%%%%%%%%%%%%%%%%%%%

  \subsection*{$L$-convex polyominoes. }
\begin{figure}[h]
    \centering
    \subfloat[]{\label{fig:L convex}\includegraphics[scale=0.5]{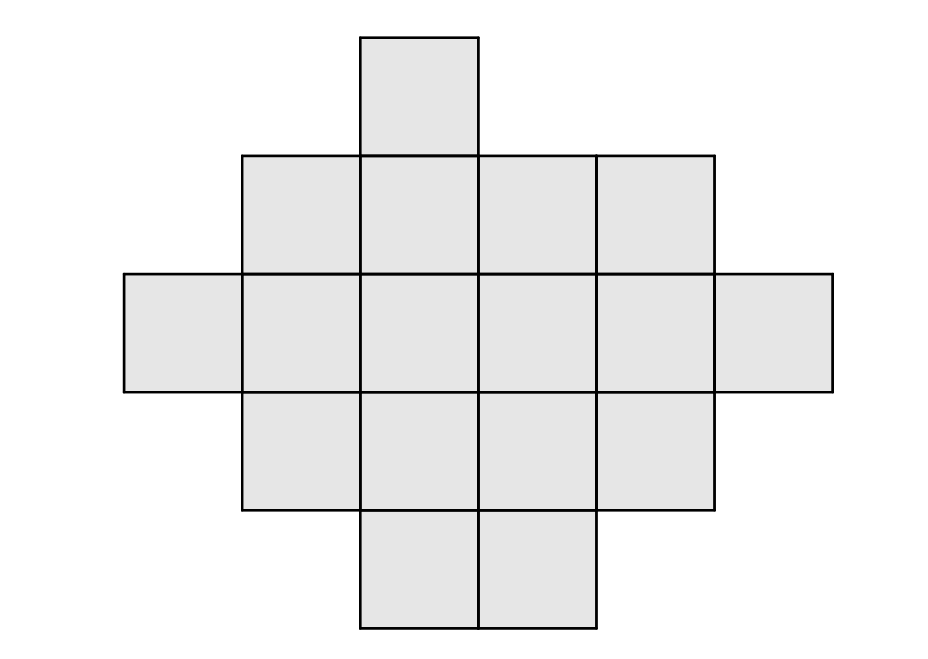}}\qquad\qquad
    \subfloat[]{\label{fig:no L-convex}\includegraphics[scale=0.5]{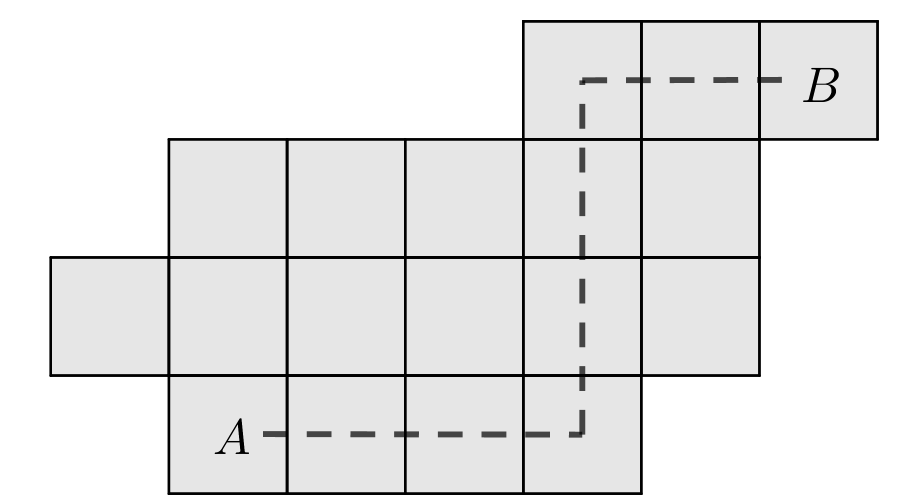}}
    \caption{An $L$-convex polyomino on the left and a $2$-convex one on the right}
\end{figure}

A convex polyomino $\cP$ is called \textit{$k$-convex} if any two cells in $\cP$ can be connected by a path of cells contained in $\cP$ with at most $k$ changes of direction. A notable case occurs when $k=1$, yielding the class of \emph{$L$-convex polyominoes}. Figures \ref{fig:L convex} and \ref{fig:no L-convex} display, respectively, an $L$-convex polyomino and one that is not $L$-convex but $2$-convex.

It is observed in \cite{EneHerzogQureshiRomeo2021} that an $L$-convex polyomino,
after a suitable permutation of its rows and columns, becomes a Ferrers diagram
(see Figure~\ref{fig:Ferrer diagram and stack}(A)).  
By Section~\ref{Section:SimplePolyominoes}, the coordinate ring $K[\Pc]$ of such a
polyomino is isomorphic to the edge ring of a weakly chordal bipartite graph; in the
Ferrers case, this graph is a Ferrers graph.  
The Castelnuovo–Mumford regularity of edge rings of Ferrers graphs is computed
in \cite[Theorem~5.7]{CorsoNagel}.  
This leads to the following result.

\begin{Theorem}\cite[Theorem 3.3]{EneHerzogQureshiRomeo2021}
Let $\cP$ be an $L$-convex polyomino. Then the Castelnuovo–Mumford regularity of $K[\cP]$ coincides with the rook number of $\cP$.
\end{Theorem}

%The proof is essentially combinatorial and relies on two main ingredients: \cite[Theorem~5.7]{CorsoNagel}, which computes the regularity of the toric ideal of a Ferrers graph in terms of the associated Young diagram, and the fact that every polyomino ideal of an $L$-convex polyomino $\cP$ can be viewed as the toric ideal of a Ferrers graph after a suitable rearrangement of the cells of $\cP$ (Theorem~3.1). 

In addition, the authors also discuss the Gorenstein property:

\begin{Theorem}\cite[Theorem 4.3]{EneHerzogQureshiRomeo2021}
Let $\cP$ be an $L$-convex polyomino, and let $\cP_0, \cP_1, \ldots, \cP_t$ be the derived sequence of $L$-convex polyominoes of $\cP$. Then the following conditions are equivalent:
\begin{enumerate}[(a)]
    \item $K[\cP]$ is Gorenstein.
    \item For $0 \le k \le t$, the bounding box of $\cP_k$ is a square.
\end{enumerate}
\end{Theorem}

Here, the fact that $\cP_0,\dots,\cP_t$ form a derived sequence means that $\cP_k$ is obtained from $\cP_{k-1}$ by removing a suitable maximal rectangle of $\cP_{k-1}$ and gluing the two remaining parts together so as to obtain the $L$-convex polyomino $\cP_k$ (see page 12 in \cite{EneHerzogQureshiRomeo2021}). This result generalizes the characterization of Gorenstein stack polyominoes established in \cite[Corollary~28]{Andrei} and \cite[Corollary~4.12]{Qureshi2012}, formulated in terms of the shape of the polyomino. However, in \cite[Theorem 21]{Andrei}, all convex polyominoes whose coordinate ring is Gorenstein are completely classified in terms of the associated graph introduced in Section~\ref{Section:SimplePolyominoes}.
%%%%%%%%%%%%%%%%%%%%%%%%%%%%%%%%%%%%%%%%%%%%%%%%%%%%%%%%%%%%%%%%%%%%%%%%%%%%%%%%%%%%%%%%%%%%%%%%%%%%%%%%%%%%%%%%%%%%%%%%%%%%%%%%%%%%%%%%%%%%%%%%%%%%%%%%%%%%%%%%%%%%%%%%%%%%%%%%%%%%
\subsection*{Simple thin polyominoes}
The relationship between rook theory and the Hilbert–Poincaré series is further developed in the work of Rinaldo and Romeo \cite{RR}.  
For simple thin polyominoes, they obtain an explicit combinatorial description for the $h$-polynomial of $K[\cP]$: it is given by the rook polynomial of $\cP$.  
This provides a second broad family of polyominoes for which algebraic invariants can be read directly from rook configurations.

\begin{Theorem}\cite[Theorem 3.13, Corollary 3.14]{RR}
    Let $\cP$ be a simple thin polyomino. Then the $h$-polynomial and the regularity of $K[\cP]$ coincide with the rook polynomial and the rook number of $\cP$.
\end{Theorem}
The proof relies on the standard decomposition of Hilbert--Poincar\'e series of a homogeneous ideal $I$ of a graded ring $R$ via the the short exact sequence
\[
0 \longrightarrow R/(I : f) \xrightarrow{\cdot f} R/I \longrightarrow R/(I, f) \longrightarrow 0,
\]
where $f \in R$ is a homogeneous element of degree $d$, which yields
\[
\mathrm{HP}_{R/I}(t) = \mathrm{HP}_{R/(I,f)}(t) + t^{d}\,\mathrm{HP}_{R/(I : f)}(t).
\]

Applied to a simple thin polyomino $\cP$, this produces the recursive formula
\[
\mathrm{HP}_{K[\cP]}(t)
= \frac{1}{1 - t}
\left(
    \mathrm{HP}_{K[\cP']}(t)
    + \frac{t}{(1 - t)^{r-1}} \, \mathrm{HP}_{K[\cP'']}(t)
\right),
\]
where $\cP'$ is obtained by deleting the leaf cells of $\cP$, and $\cP''$ is formed by removing the maximal interval containing the leaf and gluing the remaining two pieces (see \cite[Definitions~3.3 and~3.4]{RR}).  
Induction on this expression shows that
\[
h_{K[\cP]}(t) = r_{\cP}(t),
\]
without the need to compute the rook polynomial or the $h$-polynomial explicitly.  
Since $K[\cP]$ is Cohen--Macaulay, this also implies
\[
\reg K[\cP] = r(\cP).
\]

These observations motivate the following conjecture.

\begin{Conjecture}[Conjecture 4.5, Question 4.6]\cite{RR}
     Let $\cP$ be a polyomino. 
     \begin{itemize}
         \item $\cP$ is thin if and only if $r_{\cP}(t) = h_{K[\cP]}(t)$.
         \item Moreover, is $\reg K[\cP] = r(\cP)$?
     \end{itemize}
\end{Conjecture} 
 
\noindent They also investigate the Gorenstein property. This property is strongly connected to the $h$-polynomial thanks to a classical result of Stanley \cite{S}: indeed, a $K$-algebra that is a domain is Gorenstein if and only if its $h$-polynomial is palindromic. This condition is then translated into a structural property of the polyomino, known as the $S$-property (see Definition~4.1 therein). 

\begin{Theorem}\cite[Theorem 4.2]{RR}
     Let $\cP$ be a simple thin polyomino. Then $K[\cP]$ is Gorenstein if and only if $\cP$ has the $S$-property.
\end{Theorem}
 
Building on this characterization, Kummini and Veer prove in \cite{KV2} that
simple thin polyominoes with the $S$-property satisfy the Charney--Davis
conjecture.  In particular,
\[
(-1)^{\left\lfloor \deg (h_{K[\cP]}(t))/2 \right\rfloor}\, h_{K[\cP]}(-1) \ge 0.
\]

The strategy developed by Rinaldo and Romeo for simple thin polyominoes extends beyond the simple case. In \cite{CNU3}, Cisto, Navarra, and Utano proved that a similar Hilbert--Poincar\'e decomposition can be carried out for prime closed path polyominoes. The case of non-prime closed paths was then completed in \cite{CNJ}, where the approach relies on combining suitable exact sequences with Gr\"obner basis descriptions of various intermediate subpolyominoes.

\begin{Theorem}
\cite[Theorem~5.5]{CNU3}, \cite[Theorem~4.18(2)]{CNJ}, 
\cite[Theorem~13]{Navarra2025}
Let $\cP$ be a closed path polyomino.  
Then the $h$-polynomial of $K[\cP]$ coincides with the rook polynomial $r_\cP(t)$, 
and $\reg K[\cP]=r(\cP)$.
\end{Theorem}

A different proof, applicable also to weakly closed paths, is given in
\cite{Navarra2025}, and a same technique also works for grid polyominoes \cite{DinuNavarra2023}.

The Gorenstein property is likewise completely characterized:

\begin{Theorem}
\cite[Theorem~5.7]{CNU3}, \cite[Theorem~4.18(3)]{CNJ},
\cite[Theorem~13]{Navarra2025}
Let $\cP$ be a closed path polyomino.  
Then $K[\cP]$ is Gorenstein if and only if all maximal blocks of $\cP$ have rank~$3$.
\end{Theorem}

More recently, Kummini and Veer \cite{KV} proved a partial converse to the
Rinaldo--Romeo conjecture for the class of convex polyominoes whose vertex set is a sublattice of $\mathbb{N}^2$:

\begin{Theorem} \cite[Theorem~1]{KV}
Let $\cP$ be a convex polyomino whose vertex set is a sublattice of $\mathbb{N}^2$.
If $r_{\cP}(t)=h_{K[\cP]}(t)$, then $\cP$ must be thin.
\end{Theorem}

The proof exploits the distributive lattice structure of such polyominoes: the Hilbert–Poincaré series can be expressed in terms of maximal chains with $k$-descents, following the framework of Björner–Garsia–Stanley \cite{BGS}. When $\cP$ contains a square tetromino, two different $2$-rook configurations need not correspond to distinct maximal chains.  For instance, if we consider the square tetromino $\cS$, then
\[
r_{\mathcal{S}}(t)=1+4t+2t^2 
\qquad\text{while}\qquad
h_{K[\mathcal{S}]}(t)=1+4t+t^2.
\]

\begin{figure}[h]
	\centering
	\includegraphics[width=0.9\columnwidth]{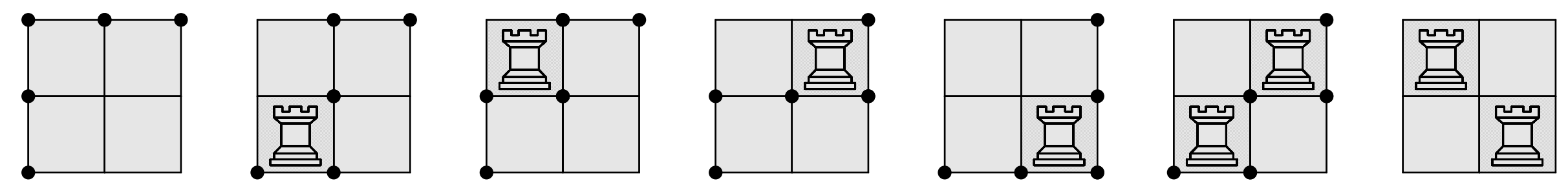}
	\caption{Rook-configurations in a square tetromino.}
    \label{Fig: equivalent tetromino}
\end{figure}

In particular, once the $1$-rook configurations and the $2$-rook configuration from the left in Figure~\ref{Fig: equivalent tetromino} are associated with maximal chains, the configuration on the right remains unassigned and cannot be matched with any maximal chain.

This observation suggests how to overcome the issue: the two $2$-rook configurations should be regarded as \textit{equivalent}, since one can be obtained from the other by moving a rook from a diagonal to an anti-diagonal position, or vice versa.

This motivates the introduction of an equivalence relation that captures this
phenomenon, indicating that in the non-thin case the appropriate combinatorial
object is not the rook polynomial itself, but a refinement of it, namely the \emph{switching rook polynomial}.  
This refined polynomial is designed precisely to account for the ambiguity in
rook configurations described above and will be introduced in the next subsection.

\subsection{Switching rook polynomial} Observe that the collection $\bigcup_{j=0}^{r(\cP)}\cR_j(\cP)$ forms a simplicial complex, known as the \textit{chessboard complex} of $\cP$. Two non-attacking rooks in $\cP$ are said to be in \textit{switching position} (or are \textit{switching rooks}) if they occupy cells that are diagonally (or anti-diagonally) opposite within an \textit{inner interval} $I$ of $\cP$, denoted $\cP_I$. In this situation, we say that the rooks are in a \textit{diagonal} (or, respectively, \textit{anti-diagonal}) position. Fix $j\in \{0,\dots, r(\cP)\}$, and let $F\in \cR_j(\cP)$. Consider two switching rooks $R_1$ and $R_2$ in $F$, positioned diagonally (or anti-diagonally) in $\cP_I$ for some inner interval $I$. Let $R_1'$ and $R_2'$ be the rooks occupying the anti-diagonal (or diagonal, respectively) cells of $\cP_I$. Then the set $(F\backslash \{R_1, R_2\}) \cup \{R_1', R_2'\}$ also lies in $\cR_j(\cP)$. This operation of replacing $R_1$ and $R_2$ with $R_1'$ and $R_2'$ is called a \textit{switch of $R_1$ and $R_2$}. 

This induces an equivalence relation $\sim$ on $\cR_j(\cP)$: we write $F_1\sim F_2$ if $F_2$ can be obtained from $F_1$ by a sequence of switches. In this case, we say that $F_1$ and $F_2$ are \textit{equivalent with respect to $\sim$} (or \textit{equal up to switches}). The following figure shows four $3$-rook configurations that are equivalent under~$\sim$.

\begin{figure}
	\hfill \begin{minipage}[t]{0.22\linewidth}
		\includegraphics[width=0.69\linewidth]{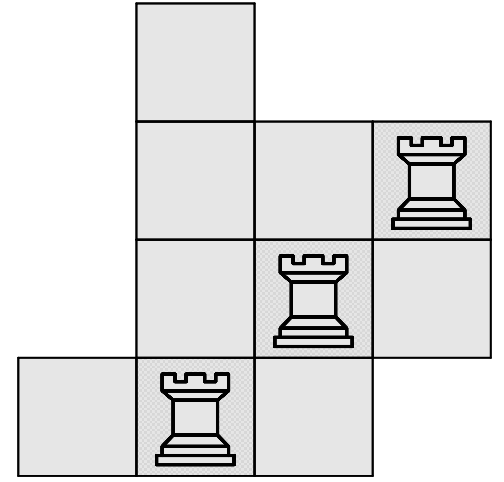}
	\end{minipage}
	\begin{minipage}[t]{0.22\linewidth}
		\includegraphics[width=0.69\linewidth]{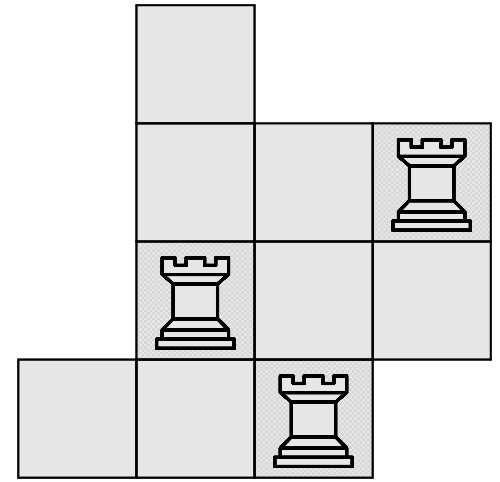}
	\end{minipage}
	\begin{minipage}[t]{0.22\linewidth}
		\includegraphics[width=0.69\linewidth]{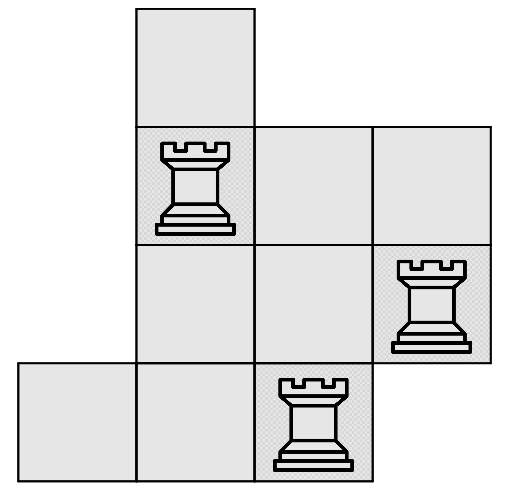}
	\end{minipage}
	\begin{minipage}[t]{0.22\linewidth}
		\includegraphics[width=0.69\linewidth]{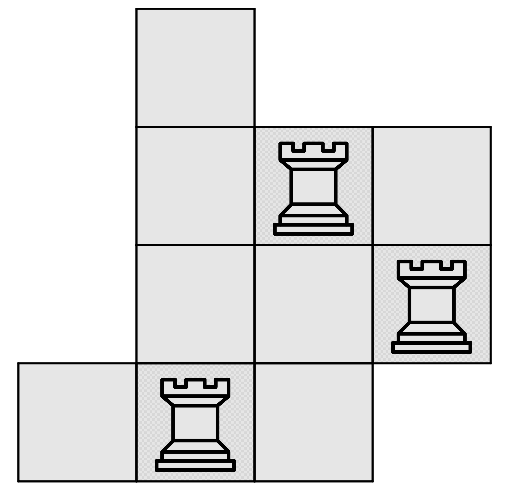}
	\end{minipage}
	\caption{Four arrangements of $3$ non-attacking rooks, equivalent under $\sim$}
\end{figure}

Let $\tilde{\cR}_j(\cP)=\cR_j(\cP)/\sim$ denote the set of equivalence classes. We set $\tilde{r}_j(\cP)=\vert \tilde{\cR}_j(\cP)\vert$ for $j\in \{0,\dots, r(\cP)\}$, with the convention that $\tilde{r}_0(\cP)=1$. The \textit{switching rook polynomial} of $\cP$ is defined as the polynomial in $\mathbb{Z}_{> 0}[t]$
$$\tilde{r}_{\cP}(t)=\sum_{j=0}^{r(\cP)}\tilde{r}_j(\cP)t^j.$$

By using different algebraic–combinatorial methods, several classes of non-thin polyominoes $\cP$ with at most one hole have been studied, and it has been shown that the $h$-polynomial of $K[\cP]$ agrees with the switching rook polynomial of $\cP$, while the regularity of $K[\cP]$ equals the rook number of $\cP$.

\subsection*{Parallelogram polyominoes and planar distributive lattices} In \cite{QureshiRinaldoRomeo2022}, Qureshi, Rinaldo, and Romeo study the Hilbert–Poincar\'e series of parallelogram polyominoes, showing that 

\begin{Theorem}\cite[Theorem 3.5, Corollary 3.13]{QureshiRinaldoRomeo2022}\label{Thm:parallelogram}
    Let $\cP$ be a parallelogram polyomino. Then the $h$-polynomial and the regularity of $K[\cP]$ coincide with the switching rook polynomial and the rook number of $\cP$.
\end{Theorem}

\noindent A parallelogram polyomino can be viewed as a planar distributive lattice and, as mentioned earlier, its Hilbert–Poincar\'e series is known and described in terms of maximal chains with $k$-descents (see \cite{BGS}). As a substantial refinement of \cite{KV}, a bijective correspondence between maximal chains with $k$-descents and arrangements of $k$ non-attacking rooks, up to switches, is established in Proposition~3.11 and Lemma~3.12, thereby proving that $\tilde{r}_\cP(t)=h_{K[\cP]}(t)$ (Theorem~3.5). Moreover, by computational methods, the latter result is proved for all simple polyominoes of rank at most $11$. 

\begin{Theorem}\cite[Theorem 3.4]{QureshiRinaldoRomeo2022}
    Let $\cP$ be a simple polyomino with rank at most $11$. Then the $h$-polynomial and the regularity of $K[\cP]$ coincide with the switching rook polynomial and the rook number of $\cP$.
\end{Theorem}

\noindent This leads to a first conjecture, namely that for any simple polyomino $\cP$, one should have $\tilde{r}_\cP(t)=h_{K[\cP]}(t)$. In \cite{NQR1}, this correspondence is reasonably extended to all collections of cells. We also remark that Theorem~\ref{Thm:parallelogram} allows one to recover the result of \cite{KV}.

\begin{Remark}\rm
Let $\cP$ be a convex polyomino whose vertex set is a sublattice of $\NN^2$ such that $r_{\cP}(t)=h_{K[\cP]}(t)$. Then, by \cite[Theorem 3.5]{QureshiRinaldoRomeo2022}, we have $\tilde{r}_{\cP}(t)=r_{\cP}(t)$, which implies that $\cP$ must be thin. 
\end{Remark}

\noindent Hibi \cite{H} characterized Gorenstein simple planar distributive lattices by proving that they are Gorenstein if and only if all maximal chains have the same length. The corresponding characterization for parallelogram polyominoes is reformulated in terms of the so-called $S$-property (see how \cite[Definition~4.1]{QureshiRinaldoRomeo2022} generalizes \cite[Definition~4.1]{RR}) in \cite[Theorem~4.10]{QureshiRinaldoRomeo2022}, where a description in terms of Motzkin paths is also provided (see \cite[Corollary~4.13]{QureshiRinaldoRomeo2022}).

\subsection*{Shellable flag simplicial complexes of polyominoes.} 
\begin{figure}[h]
	\centering
	\includegraphics[scale=0.5]{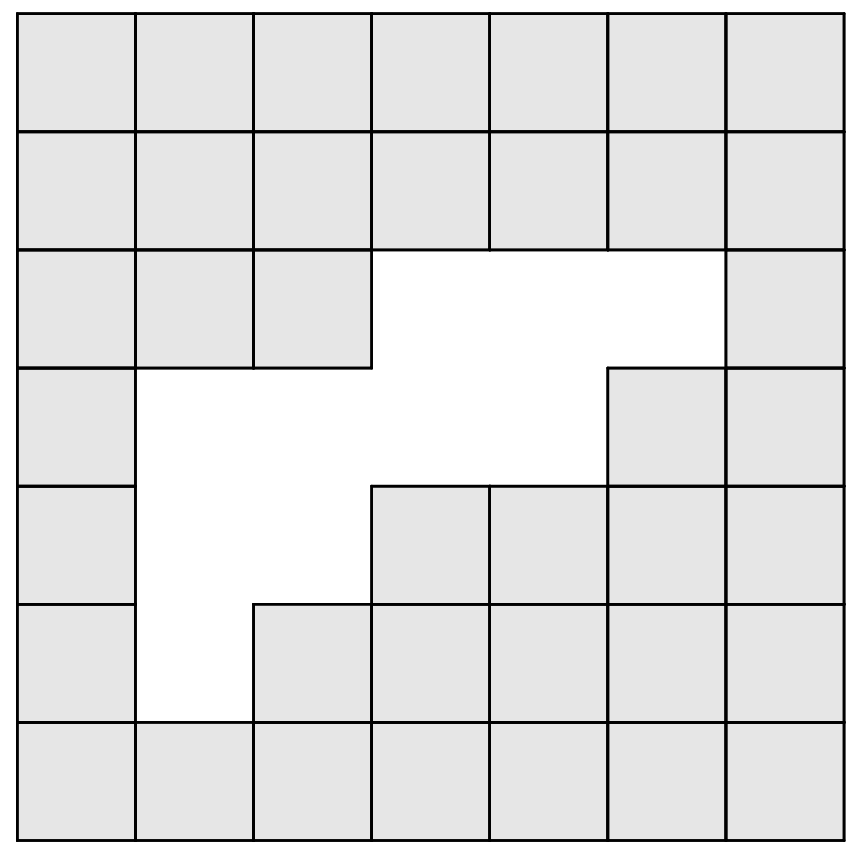}
	\caption{Frame polyomino.}
	\label{fig:frame}
\end{figure}

In \cite{JahangirNavarra2024}, \emph{frame polyominoes} are introduced by Jahangir and Navarra. A frame polyomino is a non-simple polyomino obtained by removing the cells of a parallelogram polyomino from a rectangle (see Figure~\ref{fig:frame}). 

\noindent They study the associated Hilbert–Poincar\'e series by means of a new method based on a well-known result of McMullen–Walkup concerning the $h$-vector of a shellable simplicial complex, which we recall below. 

\begin{Theorem}\cite[Corollary 5.1.14]{Bruns_Herzog}
    Let $\Delta$ be a $d$-dimensional shellable simplicial complex with shelling $F_1, F_2, \ldots, F_m$. For $j = 1, \ldots, m$, let $r_j$ be the number of facets of $\langle F_j \rangle \setminus \langle F_1, \ldots, F_{j-1} \rangle$, and set $r_0 = 0$. Then
\[
h_i = \vert\{\, j \mid r_j = i \,\} \vert \qquad \text{for } i = 0, \ldots, d.
\]
In particular, up to reordering, the numbers $r_j$ do not depend on the chosen shelling.
\end{Theorem}

\noindent The paper analyses the flag simplicial complex arising from the initial ideal of $I_{\cP}$ with respect to the reverse lexicographic order induced by the natural total order on the vertices. To achieve this, a new combinatorial notion, namely the \emph{step} of a face (see Definition~3.3 therein), is introduced, and a complete description of the facets with $k$-steps is provided in terms of the vertices of the frame polyomino (see Discussion~3.9 therein). A crucial result is then established to prove the shellability of $\Delta_{\cP}$.

\begin{Theorem}\cite[Theorem 3.12]{JahangirNavarra2024}
Let $\cP$ be a frame polyomino and let $\Delta(\cP)$ be the simplicial complex attached to $\cP$.  
Suppose that $\mathcal{F}_{\cP}$ is lexicographically ordered in descending order and denote by $<_{\mathrm{lex}}$ such a order. Consider a facet $F \neq F_0$ of $\Delta(\cP)$ and set $\cS(F) = \{\, G \in \mathcal{F}_{\cP} : F <_{\mathrm{lex}} G \,\}$ and 
\[
K_{F} = \{\, F \setminus \{v\} : v \text{ is the lower-right corner of a step of } F \,\}.
\]
Then:
\begin{enumerate}
    \item $S(F) \cap F = K_{F}$, and in particular, $\mathcal{F}_{\cP}$ forms a shelling order of $\Delta(\cP)$;
    \item the $i$-th coefficient of the $h$-polynomial of $K[\cP]$ is the number of facets of $\Delta(\cP)$ having $i$ steps.
\end{enumerate}
\end{Theorem}

\noindent Another important result is presented in \cite[Theorem 4.6]{JahangirNavarra2024}, providing a bijection between facets with $k$-steps and arrangements of $k$ non-attacking rooks, up to switches. This leads to the main theorem of the paper:

\begin{Theorem}\cite[Theorem 4.7, Corollary 4.8]{QureshiRinaldoRomeo2022}
    Let $\cP$ be a frame polyomino. Then the $h$-polynomial and the regularity of $K[\cP]$ coincide with the switching rook polynomial and the rook number of $\cP$.
\end{Theorem} 

\noindent The conjecture \cite[Conjecture~3.2]{QureshiRinaldoRomeo2022} was formulated originally for simple polyominoes. Since a frame polyomino is non-simple, it was natural to expect that the conjecture could be extended to arbitrary polyominoes, as proposed in \cite[Conjecture~4.9]{JahangirNavarra2024}.

The approach employed by Jahangir and Navarra in \cite{JahangirNavarra2024} is further extended in \cite{DinuNavarra2025} to obtain the corresponding results for grid polyominoes, and in \cite{Navarra2025} for closed path and weakly closed path polyominoes (where a suitable monomial order is chosen to guarantee the shellability of $\Delta_{\cP}$). 

\subsection*{Convex collections of cells with quadratic Gr\"obner basis.} In \cite{NQR1}, Navarra, Qureshi, and Rinaldo provide an algorithm to compute the switching rook polynomial of a collection of cells (see \cite{N1}) and, by using it, establish the following:

\begin{Theorem}\label{thm: computational thm}
Let $\cP$ be a collection of cells. Then $h_{K[\cP]}(t)$ coincides with the switching rook polynomial of $\cP$, and $\mathrm{reg}(K[\cP])$ equals the rook number of $\cP$ in the following cases:
\begin{itemize}
    \item when $\cP$ is a collection of cells of rank at most $10$;
    \item when $\cP$ is a polyomino of rank at most $12$.
\end{itemize}
\end{Theorem}

\noindent Motivated by this evidence, they conjecture that the correspondence holds in general:

\begin{Conjecture}\label{conj}
Let $\cP$ be a collection of cells. Then the switching rook polynomial of $\cP$ coincides with the $h$-polynomial of $K[\cP]$, and the rook number of $\cP$ equals the regularity of $K[\cP]$.
\end{Conjecture}

\noindent The second part of the paper proves the above conjecture for convex collections of cells whose inner $2$-minor ideal has a quadratic Gr\"obner basis with respect to $<_{\rev}$ and $<_{\lex}$, where $<_{\rev}$ and $<_{\lex}$ denote the reverse lexicographic and lexicographic orders on $S_{\Pc}$ induced by the natural total order on its variables: for $x_a, x_b \in S_{\Pc}$ with $a = (i,j)$ and $b = (k,\ell)$, we set $x_a > x_b$ if $i > k$, or if $i = k$ and $j > \ell$.

\subsection*{Domino-stability, palindromicity and Gorensteiness.} In \cite{NQR2}, the palindromicity of the switching rook polynomial is studied. The central notion introduced there is the \emph{domino stability} of a collection of cells, and a complete characterization is given of all collections of cells whose switching rook polynomial is palindromic:

\begin{Theorem}\cite[Theorem 5.1]{NQR2}
Let $\cP$ be a collection of cells and let $\tilde{r}_{\cP}(t)$ denote the switching rook polynomial of $\cP$. Then $\tilde{r}_{\cP}(t)$ is palindromic if and only if $\cP$ is domino-stable.
\end{Theorem}

\noindent Consequently, by Stanley’s classical result \cite{S}, domino stability provides a characterization of the Gorenstein property of $K[\cP]$ when $K[\cP]$ is a domain, and a necessary condition when it is not. Computational evidence obtained with \cite{N1} and \cite{N2} confirms the following:

\begin{Proposition}\cite[Proposition 5.5]{NQR2}
Let $\cP$ be a domino-stable collection of cells with rank less than or equal to $10$, or a domino-stable polyomino with rank less than or equal to $12$. Then $K[\cP]$ is Gorenstein.
\end{Proposition}

\noindent This leads to the following conjecture:

\begin{Conjecture}\cite[Conjecture 5.6]{NQR2}
Let $\cP$ be a collection of cells. Then $K[\cP]$ is Gorenstein if and only if $\cP$ is domino-stable.
\end{Conjecture}

\section{Canonical Module, Pseudo-Gorensteiness and Levelness}\label{Sectio: canonical module}

% There are two relevant generalizations of Gorenstein rings: level rings 
% and pseudo-Gorenstein rings, and both properties have been studied for polyomino ideals. 
Let $R$ be a standard graded Cohen–Macaulay $K$-algebra with canonical 
module $\omega_R$. Then $R$ is Gorenstein if and only if its canonical module is cyclic, 
and hence generated in a single degree. This condition on $\omega_R$ may be weakened in different ways, in particular:
\begin{enumerate}
    \item if one only requires that all generators of $\omega_R$ have the same degree, then $R$ is called a \emph{level} ring;
    \item if one requires that $\omega_R$ has a unique generator of minimal degree, then $R$ is called a \emph{pseudo-Gorenstein} ring.
\end{enumerate}
Pseudo-Gorenstein rings can be studied through their Hilbert–Poincar\'e series, since 
this property occurs precisely when the leading coefficient of the $h$-polynomial is equal to $1$.

\subsection*{Pseudo-Gorenstein and level paths.} 
In \cite{RinaldoRomeoSarkar2024}, Rinaldo, Romeo and Sarkar initiated the study of the 
pseudo-Gorenstein and level properties for path polyominoes. Since for simple thin polyominoes the 
$h$-polynomial coincides with the rook polynomial, when $\cP$ is a simple path this amounts 
to characterizing those having a unique $r(\cP)$-rook configuration. 

They call a \emph{stair} a ladder $\cS$ of at least two steps such that all maximal rectangles contained in $\cS$, except the first and the last one, have rank equal to $2$. 
An \emph{odd stair} is a stair containing an odd number of maximal rectangles.

\begin{Theorem}[{\cite[Theorem 4]{RinaldoRomeoSarkar2024}}]
Let $\cP$ be a path polyomino with $\{I_1, I_2, \ldots, I_s\}$ the sequence of maximal 
rectangles of $\cP$, and let $l_k = \vert I_k \vert$ for all $1 \le k \le s$. 
Then $K[\cP]$ is pseudo-Gorenstein if and only if either $\cP$ is a cell or the following 
conditions hold:
\begin{enumerate}
    \item $l_1 = l_s = 2$ and $l_k \leq 3$ for all $2 \le k \le s - 1$;
    \item $\cP$ does not contain odd stairs.
\end{enumerate}
\end{Theorem}

To characterize level simple paths, the authors study the socle of the ring $K[\cP]$ modulo suitable linear forms, together with the structure of rook configurations in $\cP$.  
By defining a \textit{bad stair} as a stair whose number of maximal rectangles is $4$, $6$, or greater than or equal to $8$, they proved the following result.

\begin{Theorem}[{\cite[Theorem 10]{RinaldoRomeoSarkar2024}}]
Let $\cP$ be a path polyomino. Then $K[\cP]$ is level if and only if $\cP$ does not contain a bad stair.
\end{Theorem}

\subsection*{Canonical module of circle closed polyominoes.}

\begin{figure}[h]
    \centering
    \includegraphics[width=0.3\linewidth]{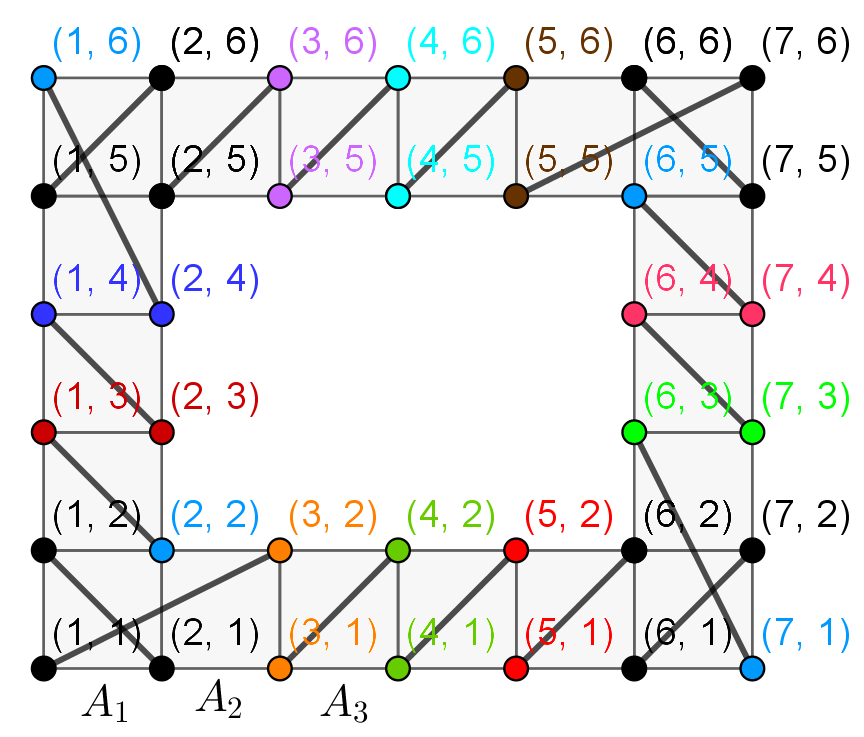}
    \caption{A circle closed path.}
    \label{fig:circle}
\end{figure}

In \cite{DinuNavarra2025}, Dinu and Navarra provide a combinatorial description of the canonical module of the coordinate ring of circle closed path polyominoes. A circle closed path is a closed path having exactly four maximal rectagles (see Figure \ref{fig:circle}, for instance). 

As a preliminary step, they establish the following result.

\begin{Theorem}[{\cite[Theorem 4.8]{DinuNavarra2025}}]
Let $\cP$ be a closed path polyomino. Then $I_{\cP}$ is of Kőnig type.
\end{Theorem}

Based on this, one can define two ideals: a binomial ideal $J(\cP)$, arising from the Kőnig-type property, and a monomial ideal $K(\cP)$, constructed from suitable vertices of $\cP$. In Figure~\ref{fig:circle}, one can observe the following:
(1) the vertices involved in the monomial generators of $K(\cP)$—in particular, the variables appearing in a generator correspond to the blue vertices and, for each pair of vertices of the same color, the generator contains exactly one of them, subject to the restriction that no two chosen vertices may lie in a diagonal position;
(2) the leading terms, represented by the “diagonal lines”, of the binomials of $I_{\cP}$, which are selected as generators of $J(\cP)$.\\
Using a classical result from linkage theory, they prove the following.

\begin{Theorem}[{\cite[Theorem 5.4, Corollary 5.17]{DinuNavarra2025}}]
Let $\cP$ be a circle closed path. Denote by $\omega_{K[\cP]}$ the canonical module of $K[\cP]$. Then
$$\omega_{K[\cP]} \cong \frac{(J(\cP)+K(\cP))}{J(\cP)}$$
where $J(\cP)$ and $K(\cP)$ are the ideals defined in \cite[Definition 5.2]{DinuNavarra2025}.\\
Moreover, $K[\cP]$ is a level ring.
\end{Theorem}

\subsection*{Fuss--Catalan numbers as the Cohen--Macaulay type of certain Ferrers diagrams}

For $p \geq 1$, let $u_1, \ldots, u_p$ and $r_1, \ldots, r_p$ be integers. We denote by 
\[
\cF=\cP\begin{pmatrix}
u_1 & u_2 & \cdots & u_p \\
r_1 & r_2 & \cdots & r_p
\end{pmatrix}
\]
the Ferrers diagram whose first $r_1$ columns consist of $u_1$ cells, the next $r_2$ columns consist of $(u_1+u_2)$ cells, and, in general, the next $r_k$ columns consist of $\sum_{i=1}^k u_i$ cells, up to the final $r_p$ columns of $\sum_{i=1}^p u_i$ cells.  

\begin{figure}[h]
    \centering
    \includegraphics[width=0.3\linewidth]{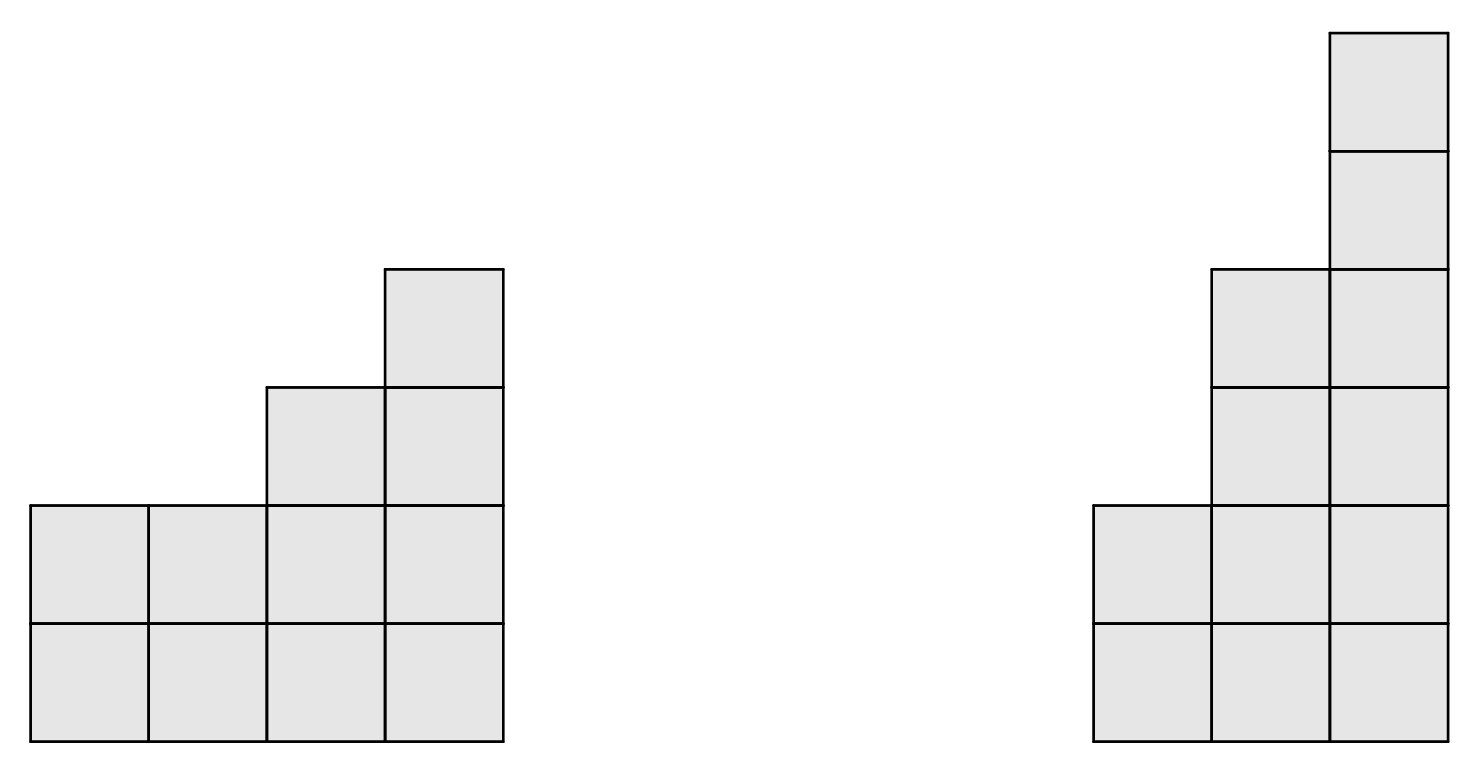}
    \caption{Ferrers diagrams.}
    \label{fig:ferrer for Fuss}
\end{figure}

For instance, Figure~\ref{fig:ferrer for Fuss} displays $\cP\begin{pmatrix}
2 & 1 & 1 \\
2 & 1 & 1 
\end{pmatrix}
$ on the left  and $\cP\begin{pmatrix}
2 & 2 & 2 \\
1 & 1 & 1
\end{pmatrix}
$ on the right.

Ştefan provided an explicit description of the Cohen--Macaulay type for a particular class of Ferrers diagrams, expressed in terms of the Fuss--Catalan numbers
\[
C_p(n) = \frac{1}{(n-1)p+1} \binom{np}{p}.
\]

\begin{Theorem}[{\cite[Theorem 10, Corollary 11]{AlinStefan}}]
Let
\[
\cP_1 = \cP\begin{pmatrix}
u_1 & u_2 & \cdots & u_p \\
r_1 & r_2 & \cdots & r_p
\end{pmatrix}
\qquad\text{and}\qquad
\cP_2 = \cP\begin{pmatrix}
u_1 & u_2 & \cdots & u_p \\
s_1 & s_2 & \cdots & s_p
\end{pmatrix},
\]
where $u_1=\cdots=u_p=n$, $r_1=\cdots=r_p=t$, and $s_1=\cdots=s_p=n-t$.  
Then the two non-isomorphic coordinate rings $K[\cP_1]$ and $K[\cP_2]$ have the same Cohen--Macaulay type. Moreover, for $t=1$ this type coincides with the Fuss--Catalan number $C_{p+1}(n)$.
\end{Theorem}

\noindent
The proof relies on the classical theorem of Danilov--Stanley, which describes the canonical module $\omega_R$ of a semigroup ring $R=K[\cA]$, where $\cA$ is a collection of lattice points in $\mathbb{Z}^n$, in terms of the polyhedral cone generated by~$\cA$.

\section{A package for \texttt{Macaulay2}: \texttt{PolyominoIdeals}}\label{Section: Package}

In \cite{Package_M2}, Cisto, Jahangir, and Navarra developed a package for the computer algebra system \texttt{Macaulay2} \cite{M2} designed to work with collections of cells and their associated binomial ideals.

Since a cell in the lattice $\mathbb{Z}^2$ is uniquely determined by its lower-left corner, it is natural to encode a collection of cells as a list of such corners. For instance, a square tetromino may be encoded as  
\texttt{Q = \{\{1, 1\}, \{1, 2\}, \{2, 1\}, \{2, 2\}\}}. The command \texttt{cellCollection Q} then allows to create an object of type \texttt{CollectionOfCells}, that can use in the provided functions.

Below we provide the complete list of functions implemented in the package, organised according to their structural–combinatorial, rook-theoretic, and algebraic features.

\subsection{Structural–combinatorial functions}

\begin{description}

  \item[\href{https://www.macaulay2.com/doc/Macaulay2/share/doc/Macaulay2/PolyominoIdeals/html/_cell__Collection.html}{\texttt{cellCollection}}]
        Create a collection of cells.

  \item[\href{https://www.macaulay2.com/doc/Macaulay2/share/doc/Macaulay2/PolyominoIdeals/html/_polyo__Vertices.html}{\texttt{polyoVertices}}]
       Provide the vertices of a collection of cells.

  \item[\href{https://www.macaulay2.com/doc/Macaulay2/share/doc/Macaulay2/PolyominoIdeals/html/_inner__Interval.html}{\texttt{innerInterval}}]
        Check if an interval is an inner interval of a collection of cells.

  \item[\href{https://www.macaulay2.com/doc/Macaulay2/share/doc/Macaulay2/PolyominoIdeals/html/_cell__Graph.html}{\texttt{cellGraph}}]
        Provide the graph $G$ associated with a collection of cells $\{C_1,\dots,C_n\}$, where $V(G)=[n]$ and $E(G)=\{ \{i,j\}: C_i \text{ shares an edge with } C_j\}$. 

  \item[\href{https://www.macaulay2.com/doc/Macaulay2/share/doc/Macaulay2/PolyominoIdeals/html/_collection__Is__Connected.html}{\texttt{collectionIsConnected}}]
        Check whether a collection of cells is connected.

  \item[\href{https://www.macaulay2.com/doc/Macaulay2/share/doc/Macaulay2/PolyominoIdeals/html/_connected__Components__Cells.html}{\texttt{connectedComponentsCells}}]
        Provide the connected components of a collection of cells.

  \item[\href{https://www.macaulay2.com/doc/Macaulay2/share/doc/Macaulay2/PolyominoIdeals/html/_is__Row__Convex.html}{\texttt{isRowConvex}}]
        Check the row convexity of a collection of cells.

  \item[\href{https://www.macaulay2.com/doc/Macaulay2/share/doc/Macaulay2/PolyominoIdeals/html/_is__Column__Convex.html}{\texttt{isColumnConvex}}]
        Check the column convexity of a collection of cells.

  \item[\href{https://www.macaulay2.com/doc/Macaulay2/share/doc/Macaulay2/PolyominoIdeals/html/_is__Convex.html}{\texttt{isConvex}}]
        Check the convexity of a collection of cells.

  \item[\href{https://www.macaulay2.com/doc/Macaulay2/share/doc/Macaulay2/PolyominoIdeals/html/_collection__Is__Simple.html}{\texttt{collectionIsSimple}}]
        Check if a collection of cells is simple.

  \item[\href{https://www.macaulay2.com/doc/Macaulay2/share/doc/Macaulay2/PolyominoIdeals/html/_rank__Collection.html}{\texttt{rankCollection}}]
        Give the rank of a collection of cells.

  \item[\href{https://www.macaulay2.com/doc/Macaulay2/share/doc/Macaulay2/PolyominoIdeals/html/_random__Collection__With__Fixed__Rank.html}{\texttt{randomCollectionWithFixedRank}}]
        Provide a random collection of cells with fixed rank.

  \item[\href{https://www.macaulay2.com/doc/Macaulay2/share/doc/Macaulay2/PolyominoIdeals/html/_random__Collection__Of__Cells.html}{\texttt{randomCollectionOfCells}}]
        Provide a random collection of cells up to a given size.

  \item[\href{https://www.macaulay2.com/doc/Macaulay2/share/doc/Macaulay2/PolyominoIdeals/html/_random__Polyomino__With__Fixed__Rank.html}{\texttt{randomPolyominoWithFixedRank}}]
        Provide a random polyomino with fixed rank.

  \item[\href{https://www.macaulay2.com/doc/Macaulay2/share/doc/Macaulay2/PolyominoIdeals/html/_random__Polyomino.html}{\texttt{randomPolyomino}}]
        Provide a random polyomino of random size up to a prescribed bound.

\end{description}

\subsection{Rook theory functions}

\begin{description}

  \item[\href{https://www.macaulay2.com/doc/Macaulay2/share/doc/Macaulay2/PolyominoIdeals/html/_is__Non__Attacking__Rooks.html}{\texttt{isNonAttackingRooks}}]
        Check whether a rook configuration is non-attacking.

  \item[\href{https://www.macaulay2.com/doc/Macaulay2/share/doc/Macaulay2/PolyominoIdeals/html/_all__Non__Attacking__Rook__Configurations.html}{\texttt{allNonAttackingRookConfigurations}}]
        Provide the list of all non-attacking rook configurations in a collection of cells.

  \item[\href{https://www.macaulay2.com/doc/Macaulay2/share/doc/Macaulay2/PolyominoIdeals/html/_rook__Polynomial.html}{\texttt{rookPolynomial}}]
        Compute the rook polynomial of a collection of cells.

  \item[\href{https://www.macaulay2.com/doc/Macaulay2/share/doc/Macaulay2/PolyominoIdeals/html/_rook__Number.html}{\texttt{rookNumber}}]
        Compute the rook number of a collection of cells.

  \item[\href{https://www.macaulay2.com/doc/Macaulay2/share/doc/Macaulay2/PolyominoIdeals/html/_equivalence__Classes__Switching__Rook.html}{\texttt{equivalenceClassesSwitchingRook}}]
        Provide the list of equivalence classes of non-attacking rook configurations under switching.

  \item[\href{https://www.macaulay2.com/doc/Macaulay2/share/doc/Macaulay2/PolyominoIdeals/html/_switching__Rook__Polynomial.html}{\texttt{switchingRookPolynomial}}]
        Compute the switching rook polynomial of a collection of cells.

  \item[\href{https://www.macaulay2.com/doc/Macaulay2/share/doc/Macaulay2/PolyominoIdeals/html/_standard__Rook__Number.html}{\texttt{standardRookNumber}}]
        Compute the standard rook number of a collection of cells.

  \item[\href{https://www.macaulay2.com/doc/Macaulay2/share/doc/Macaulay2/PolyominoIdeals/html/_standard__Non__Attacking__Rook__Configurations.html}{\texttt{standardNonAttackingRookConfigurations}}]
        Give the list of the standard non-attacking rook configurations.

  \item[\href{https://www.macaulay2.com/doc/Macaulay2/share/doc/Macaulay2/PolyominoIdeals/html/_standard__Rook__Polynomial.html}{\texttt{standardRookPolynomial}}]
        Compute the standard rook polynomial of a collection of cells.

\end{description}

\subsection{Algebraic functions}

\begin{description}

  \item[\href{https://www.macaulay2.com/doc/Macaulay2/share/doc/Macaulay2/PolyominoIdeals/html/_polyo__Ideal.html}{\texttt{polyoIdeal}}]
        Give the ideal of inner 2-minors of a collection of cells.

  \item[\href{https://www.macaulay2.com/doc/Macaulay2/share/doc/Macaulay2/PolyominoIdeals/html/_polyo__Toric.html}{\texttt{polyoToric}}]
        Compute the toric ideal of a collection of cells.

  \item[\href{https://www.macaulay2.com/doc/Macaulay2/share/doc/Macaulay2/PolyominoIdeals/html/_polyo__Lattice.html}{\texttt{polyoLattice}}]
        Compute the lattice ideal associated with a collection of cells.

  \item[\href{https://www.macaulay2.com/doc/Macaulay2/share/doc/Macaulay2/PolyominoIdeals/html/_adjacent2__Minor__Ideal.html}{\texttt{adjacent2MinorIdeal}}]
        Provide the ideal generated by adjacent 2-minors.

  \item[\href{https://www.macaulay2.com/doc/Macaulay2/share/doc/Macaulay2/PolyominoIdeals/html/_is__Palindromic.html}{\texttt{isPalindromic}}]
        Check whether a polynomial is palindromic.

  \item[\href{https://www.macaulay2.com/doc/Macaulay2/share/doc/Macaulay2/PolyominoIdeals/html/_polyo__Matrix.html}{\texttt{polyoMatrix}}]
        Give the matrix associated with a collection of cells.

  \item[\href{https://www.macaulay2.com/doc/Macaulay2/share/doc/Macaulay2/PolyominoIdeals/html/_polyo__Matrix__Reduced.html}{\texttt{polyoMatrixReduced}}]
        Compute the reduced form of the polyomino matrix, in according to \cite{OhsugiHibiKoszul1999}.

\end{description}

%%%%%%%%%%%%%%%%%%%%%%%%%%%%%%%%%%%%%%%%%%%%%%%%%%%%%%%%%%%%%%%%%%%%%%%%%%%%%%%%%%%%%%%%%%%%%%%%%%%%%%%%%%%%%%%%%%%%%%%%%%%%%%%%%%%%%%%%%%%%%%%%%%%%%%%%%%%%%%%%%%%%%%%%%%%%%%%%%%%%%%%%%%%%%%%%%%%%%%%%%%%%%%%%%%%%%%%%%%%%%%%%%%%%%%%%%%%%%%%%%%%%%%%%%%%%%%%%%%%%%%%%%%%%%%%%%%%%%%%%%%%%%%%%

% \bibitem{HamonanganMuchtadi2024} Y.~Y.~Hamonangan and I.~Muchtadi-Alamsyah, \textit{Binomials arising from Buchberger algorithm on polyomino ideals}, Eur. J. Pure Appl. Math. \textbf{17} (2024), no.~4, 2621--2650.

% \bibitem{HHM}
% J.~Herzog, T.~Hibi, and S.~Moradi,
% \textit{Binomial ideals attached to finite collections of cells},
% Commun. Algebra (2024).

% \bibitem{HH_Book}
% J.~Herzog and T.~Hibi,
% \textit{Monomial Ideals},
% Grad. Texts in Math., vol.~260, Springer, Berlin, 2010.

% \bibitem{HO1}H. Ohsugi, T. Hibi, Special simplices and Gorenstein toric rings,J. Combin. TheorySer. A113:718–725, 2006.

% \bibitem{Sturmfels1995}
% B.~Sturmfels,
% \textit{Gröbner Bases and Convex Polytopes},
% Univ. Lecture Series, vol.~8, Amer. Math. Soc., Providence, RI, 1995.

% \bibitem{Romeo}
% F.~Romeo,
% \textit{The Stanley–Reisner ideal of the rook complex of polyominoes},
% J. Algebra Appl. (to appear), 2650003.

% \bibitem{HMP}
% J.~Herzog, F.~Mohammadi, and J.~Page,
% \textit{Measuring the non-Gorenstein locus of Hibi rings and normal affine semigroup rings},
% J. Algebra \textbf{540} (2019), 78--99.

%%%%%%%%%%%%%%%%%%%%%%%%%%%%%%%%%%%%%%%%%%%%%%%%%%%%%%%%%%%%%%%%%%%%%%%%%%%%%%%%%%%%%%%%%%%%%%%%%%%%%%%%%%%%%%%%%%%%%%%%%%%%%%
\end{document}